%% file: GA3mex__revised_.tex
\documentclass[10]{article}
\input{macros}

\usepackage{array}
\usepackage{tabularx}
\usepackage{amsfonts}
\usepackage{amssymb}
\usepackage{tikz-cd}
\usetikzlibrary{matrix,arrows,decorations.pathmorphing}
\usepackage{graphicx}   
\def\bpm{\begin{pmatrix}}
\def\epm{\end{pmatrix}}
\makeindex

\title{From Vectors to Geometric Algebra}
\author{ Sergio Ramos Ramirez,
	\\sergio1.ramos@vw.com.mx
	\\Jos\'e Alfonso Ju\'arez Gonz\'alez,
	\\alfonso.juarez@vw.com.mx
	\\Volkswagen de  M\'exico
\\  72700 San Lorenzo Almecatla,Cuautlancingo, Pue., M\'exico	
	\\Garret Sobczyk 
	\\garret\textunderscore sobczyk@yahoo.com
\\ Universidad de las Am\'ericas-Puebla
 \\ Departamento de F\'isico-Matem\'aticas
\\72820 Puebla, Pue., M\'exico
}

\begin{document}

\maketitle

\begin{abstract} Geometric algebra is the natural outgrowth of the concept of
	a vector and the addition of vectors. After reviewing the properties of the addition of vectors, a multiplication of vectors is introduced in such a way that it encodes the famous Pythagorean theorem. Synthetic proofs of theorems in Euclidean geometry can then be replaced by powerful algebraic proofs. Whereas we largely limit our attention to 2 and 3 dimensions, geometric algebra is applicable in any number of dimensions, and in both Euclidean and non-Euclidean geometries.	
 
%\smallskip
%\no {\em AMS Subject Classification:} 15A63, 15A66, 81R05, 81R25
%\smallskip

%\no {\em Keywords:}
\end{abstract}

 \section*{0\quad Introduction}

 The evolution of the concept of number, which is at the heart of mathematics,
 has a long and fascinating history that spans many centuries and the rise and
 fall of many civilizations \cite{TD1967}. 
 Regarding the introduction of negative and complex numbers, Gauss remarked in 1831, that 
 ``... these advances, however, have always been made at first with timorous and hesitating steps".
   % cite pag 189-190 number the language of science
 %As Henri Poincar\'e so eloquently put it, {\it ``Though the source be obscure, still %the stream flows on."}   \cite{TD1930} ...
In this work, we lay down for the uninitiated reader the most basic ideas and methods of geometric algebra. Geometric algebra, the natural generalization of the real and complex number 
systems to include new quantities called {\it directed numbers}, was discovered by William Kingdon Clifford (1845-1879) shortly before his
death \cite{WKC1882}. 

In Section 1, we extend the real number system $\R $ to include {\it vectors} which are {\it directed line segments} having both {\it length} and {\it direction}. Since the geometric significance of the addition of vectors, and the multiplication of vectors by real numbers or {\it scalars}, are well understood, we only provide a short review. We wish to emphasize that the concept of a vector as a directed line segment in a flat space is independent of any coordinate system, or the dimension of the space. What is important is that the {\it location} of the directed line segment in flat space is unimportant, since a vector at a point can be translated to a parallel vector at any other point, and have the same length and direction. 
 
 Section 2 deals with the {\it geometric multiplication} of vectors. 
 Since we can both {\it add} and {\it multiply} real numbers, if the real number system is to be truly extended to include vectors, then we must be able to {\it multiply} as well as to {\it add} vectors. For guidance on how to geometrically multiply vectors, we recall the two millennium  old Pythagorean Theorem relating the sides of a right angle. By only giving up the law of {\it universal commutativity} of multiplication, we discover that the product of orthogonal vectors is anti-commutative and defines a new directed number called a {\it bivector}.
   The {\it inner} and {\it outer products} are defined in terms of the {\it symmetric} and {\it anti-symmetric} parts of the geometric product of vectors, and various important relationships between these three products are investigated. 
   
 In Section 3, we restrict ourselves to the most
 basic geometric algebras $\G_2$ of the 
Euclidean plane $\R^2$, and the geometric algebra $\G_3$ of Euclidean space $\R^3$. These geometric algebras offer concrete examples and calculations based upon the familiar rectangular coordinate systems of two and three dimensional space, although the much more general discussion of the previous sections should not be forgotten. At the turn of the 19th Century, the great {\it quaternion} verses standard {\it Gibbs-Heaviside} vector algebra was fought \cite{MJC1985}. We show how the standard {\it cross product} of two vectors is the natural {\it dual} to the outer product of those vectors, as well as the relationship to other well known identities in standard vector analysis. These ideas can easily be generalized to higher dimensional geometric algebras of both Euclidean and non-Euclidean spaces, used extensively in Einstein's famous theories of relativity \cite{DL07}, and across the mathematics \cite{H/S,LP97}, and the engineering fields \cite{ECGS01,HD02}. 

In Section 4, we treat elementary ideas from analytic geometry, including the vector equation of a line and the vector equation of a plane. Along the way, formulas for the decomposition of a vector into parallel and perpendicular components to a line and plane are derived, as well as formulas for the reflection and rotation of a vector in 2, 3 and higher dimensional spaces.

In Section 5, the flexibility and power of geometric algebra is fully revealed by discussing  {\it stereographic projection} of the unit 2-sphere centered at the origin onto the Euclidean 2-plane. Stereographic projection, and its generalization to higher dimensions, has profound applications in many areas of mathematics and physics. For example, the fundamental $2$-component spinors used in quantum mechanics have a direct interpretation in the stereographic projection of the $2$-sphere \cite{Shopf2015}.

It is remarkable that almost 140 years after its discovery, this powerful geometric number system, the natural completion of the real number system to include the
concept of direction, is not universally known by the wider scientific community, although there have been many developments and applications of the language at the advanced levels in mathematics,
theoretical physics, and more recently in the computer science and robotics communities. We feel that the main reason for this regrettable state of affairs, has been the lack of a concise, yet rigorous introduction
at the most fundamental level. For this reason we pay careful attention to
introducing the inner and outer products, and developing the basic identities,
in a clear and direct manner, and in such a way that generalization to
higher dimensional Euclidean and non-Euclidean geometric algebras presents no new obstacles for the reader. We give careful references to more advanced material, which
the interested reader can pursue at their leisure.  
 	
\section{Geometric addition of vectors}
\begin{figure}[b]
	\begin{center}
		\includegraphics[scale=.08]{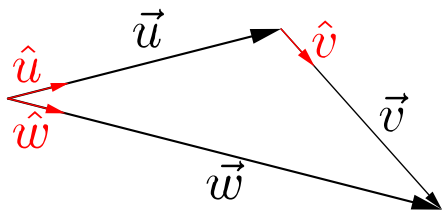}
		\caption{Vector addition.}
		\label{vectoradd}
	\end{center}
\end{figure}

	Natural numbers, or counting numbers, are used to express quantities of objects, such as 3 cows, 4 pounds, or 5 steps to north. Historically, natural numbers have been gradually extended to include fractions, negative numbers, and all numbers on the one-dimensional number line. 
	%, pictured in Figure \ref{realline}.  
{\it Vectors}, or {\it directed line-segments}, are a new kind of number which include the notion of direction. A vector $\bv=|\bv|\hat \bv$ has {\it length} $|\bv|$ and a {\it unit direction} $\hat \bv$, pictured in Figure \ref{vectoradd}. Also pictured is the sum of vectors $\bw =\bu + \bv  $.

 \begin{figure}[h]
	\centering
	\includegraphics[width=0.7\linewidth]{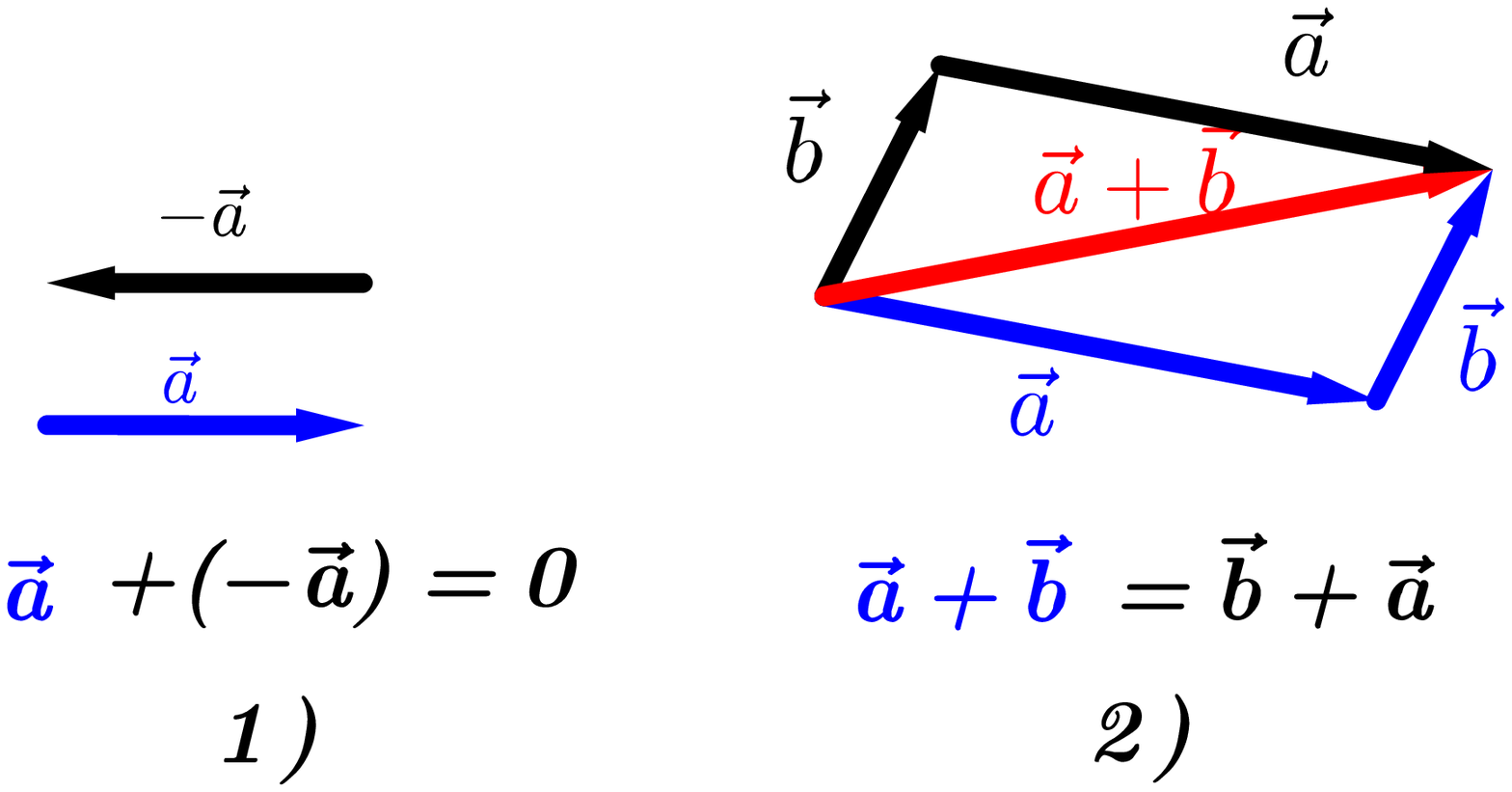}
	    	\label{prop1and2}
\end{figure}

    \begin{figure}[h]
	\centering
	\includegraphics[width=0.7\linewidth]{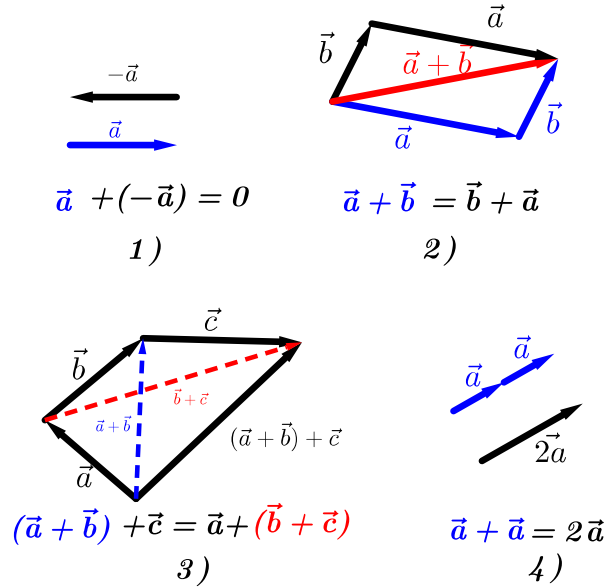}
	\caption{Geometric properties of addition of vectors.}% \textcolor{red}{In 1) replace $\ba + (-\ba)=0$ with $\bv +(-\bv )=0$.} }
	\label{prop3and4}
\end{figure}
Let $\ba $, $\bb $ and $\bc $ be vectors.
Each of the pictures in Figure \ref{prop3and4} expresses a basic geometric property of the addition of vectors, together with its translation  into a corresponding algebraic rule. For example, the {\it negative} of a vector $\ba $ is the vector $-\ba$, which has the same length as the vector $\ba$ but the {\it opposite direction} or {\it orientation}, shown in Figure \ref{prop3and4}:\ 1).  We now summarize the algebraic rules for the geometric additions of vectors, and multiplication by real numbers.

    \renewcommand{\labelenumi}{(A\arabic{enumi})}
   \begin{enumerate}
   	\item 
   	$\ba + (-\ba)= 0 \ba =  \ba 0 =0$  
   	  \hfill {\it Additive inverse of a vector}
   	
   	\item 
   	$\ba + \bb  =  \bb + \ba$  \hfill {\it Commutative law of vector addition}
   	
   	\item 
   	$(\ba + \bb) + \bc  = \ba + (\bb + \bc)  := \ba + \bb + \bc $ \hfill {\it Associative law of} 
   	\linebreak  \strut\hfill {\it vector addition}
   	
   	\item For each $\alpha \in \R $,\ 
   	$\alpha \ba =  \ba \alpha$ \hfill {\it Real numbers commute with vectors}
   	
   	\item 
   	$\ba -\bb := \; \ba + (-\bb)$ \hfill {\it Definition of vector subtraction}
   \end{enumerate}
     
     In Property (A1), the same symbol ${0}$ represents both the zero vector and the zero scalar. Property (A4), tells us that the multiplication of a vector with a real number is a commutative operation.
     Note that rules for the addition of vectors are the same as for the addition of real numbers.    
   Whereas vectors are usually introduced in
   terms of a coordinate system, we wish to emphasize that their geometric properties are independent of any coordinate system. In Section 4, we carry out explicit calculations in the geometric algebras $\G_2$ and $\G_3$, by using the usual orthonormal coordinate systems of $\R^2$ and $\R^3$, respectively.    
   
\section{Geometric multiplication of vectors} 

The geometric significance of the addition of vectors is pictured in Figures \ref{vectoradd} and \ref{prop3and4}, and formalized in the rules (A1) - (A5). But
what about the {\it multiplication} of vectors? We both add and multiply real numbers, so why can't we do the same for vectors? Let's see if we can discover how to multiply vectors in a geometrically meaningful way.  

First recall that any vector $\ba = |\ba | \hat \ba $. 
Squaring this vector, gives
\beq \ba^2 = (|\ba | \hat \ba)( |\ba | \hat \ba) = | \ba |^2 \hat \ba^2 =  | \ba |^2 . 
 \label{vecsquared} \eeq
In the last step, we have introduced the {\it new rule} that a unit vector squares to $+1$. This is always true for unit {\it Euclidean vectors}, the vectors which we are
most familiar.\footnote{{\it Space-time vectors} 
in Einstein's {\it relativity theory}, as well as vectors in other 
{\it non-Euclidean geometries}, have unit vectors with square $-1$.} With
this assumption it directly follows that a Euclidean vector squared is its {\it magnitude} or {\it length} squared, $\ba^2 = | \ba |^2 \ge 0$, and is equal to zero only when it has zero length.

 Dividing both sides of equation (\ref{vecsquared}) by
$|\ba |^2$, gives
\beq \frac{\ba^2}{|\ba |^2} = \ba \frac{\ba}{|\ba |^2}= \frac{\ba}{|\ba |^2}\ba   = 1, \label{vecinverse} \eeq
or
\[ \ba \ba^{-1} = \ba^{-1} \ba = 1 \] 
where 
 \beq \ba^{-1}:=\frac{1}{ \ba }=  \frac{\ba}{|\ba |^2} =  \frac{\hat\ba}{|\ba |} \label{inversevec} \eeq
is the {\it multiplicative inverse} of the vector $\ba $. Of course, the inverse of a vector is
only defined for nonzero vectors.
\begin{figure}[h]
	\centering
	\includegraphics[width=0.7\linewidth]{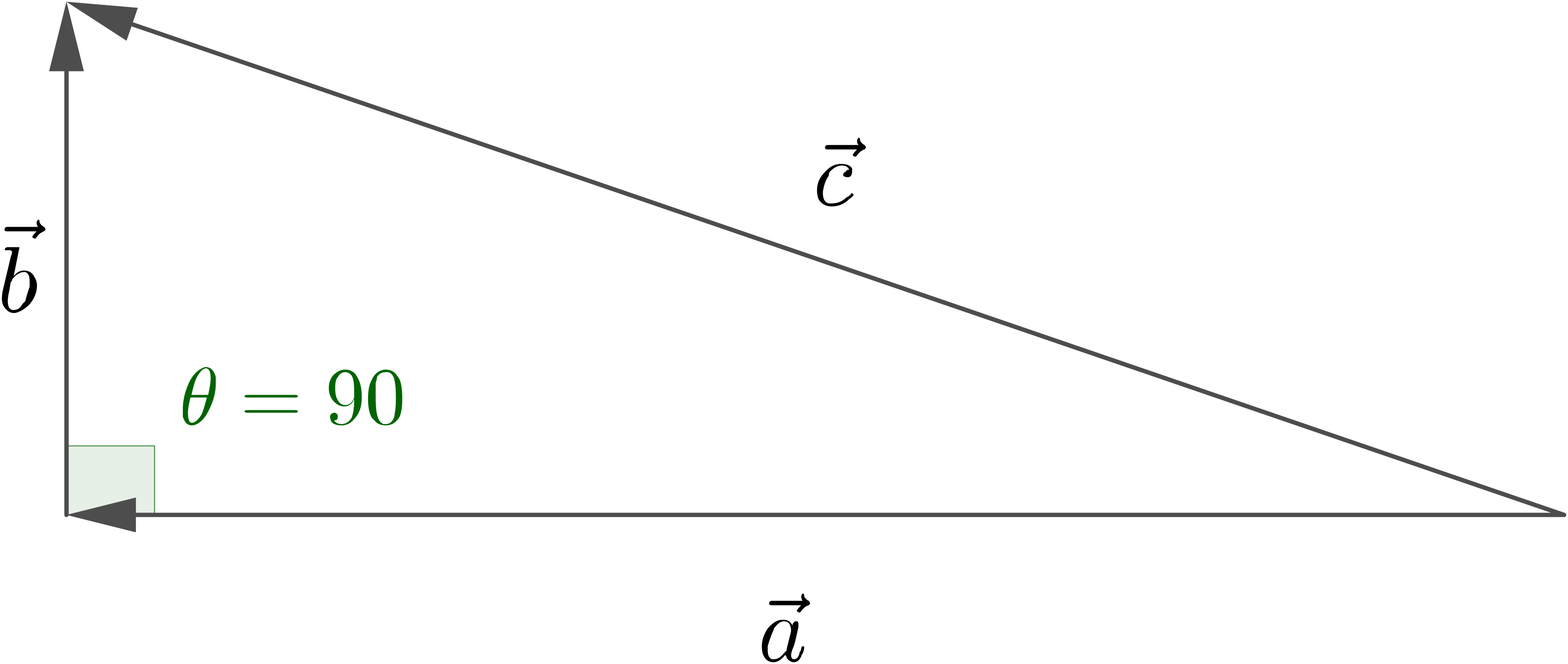}
	\caption{Right triangle with sides $\ba + \bb = \bc $. }
	\label{righttriangle}
\end{figure}

 Now consider the right triangle in Figure \ref{righttriangle}. The vectors $\ba , \bb , \bc $ along its sides satisfy the equation
 \beq \ba + \bb = \bc . \label{rtriangle} \eeq  
The most famous theorem of ancient Greek mathematics, the Pythagorean Theorem, tells
us that the lengths $|\ba|,|\bb|,|\bc|$ of the sides of this right triangle satisfy the famous relationship $|\ba|^2 + |\bb|^2 = |\bc|^2$.
Assuming the usual rules for the addition and multiplication of real numbers, except for the commutative law of multiplication, we square both sides of the vector equation (\ref{rtriangle}), to get
\[  (\ba + \bb )^2=\ba^2 +\ba \bb + \bb \ba + \bb^2=\bc^2 \ \ \iff \ \ |\ba|^2 + \ba \bb + \bb \ba + |\bb |^2= |\bc |^2,   \]
from which it follows that $\ba \bb = -\bb \ba $, if the Pythagorean Theorem is to remain valid. We have discovered that the geometric product of the orthogonal vectors $\ba $ and $\bb $ must {\it anti-commute}
if this venerable theorem is to remain true. 

\begin{figure}[h]
	\centering
	\includegraphics[width=0.7\linewidth]{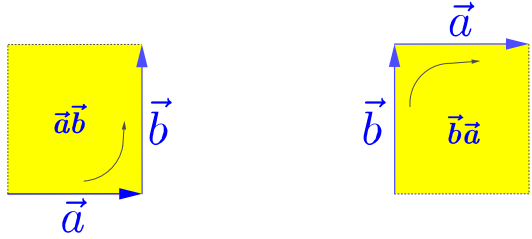}
	\caption{The bivectors $\ba \bb $ and $\bb \ba $ defined by the orthogonal vectors $\ba $ and $\bb $.}
	\label{fig:bivector}
\end{figure}

For the orthogonal vectors $\ba $ and $\bb $, let us go further and give the new quantity $\bB := \ba \bb $ the geometric interpretation of a {\it directed plane segment}, or {\it bivector}, having the direction of the plane in which the vectors lies. The bivectors $\bB $, and its {\it additive inverse} $ \bb \ba = -\ba \bb = - \bB $,
are pictured in Figure \ref{fig:bivector}. Just as the orientation of a vector is the determined by the direction of the line segment, the {\it orientations} of the bivectors $\bB = \ba \bb $ and  $-\bB = \bb \ba $ are determined by the orientation of its sides, as shown in the Figure \ref{fig:bivector}. 

We have seen that a vector $\bv = |\bv |\hat \bv $ has the unit direction $\hat \bv $ and length $|\bv |$, and that $\bv^2 = |\bv |^2$. Squaring the
bivector $\bB = \ba \bb $ gives
\beq \bB^2 = (\ba \bb)(\ba \bb)= - \ba \bb \bb \ba 
=- \ba^2 \bb^2 =-|\ba |^2|\bb |^2= -|\bB |^2, \label{bivectorsquared}\eeq
which is the {\it negative} of the area squared of the rectangle  with the sides defined by the orthogonal vectors $\ba $ and $\bb $. It follows that 
\beq \bB =|\bB| \hat \bB, \label{defbivector}  \eeq
where $|\bB |=|\ba ||\bb |$ is the area of the directed plane segment, and its direction is the {\it unit bivector} 
$\hat \bB = \hat \ba \hat \bb $, with 
\[ \hat \bB^2 =(\hat \ba \hat \bb)(\hat \ba \hat \bb)  
 =\hat \ba (\hat \bb \hat \ba) \hat \bb =-\hat \ba^2 \hat \bb^2 =   -1 . \]  

\subsection{The inner product}

\begin{figure}[h]
  	\centering
  	\includegraphics[width=0.7\linewidth]{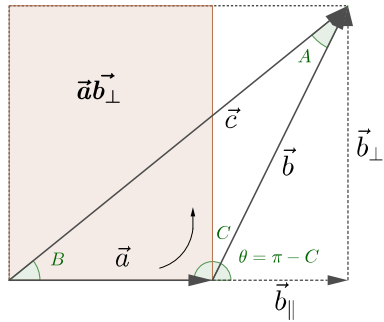}
  	\caption{Law of Cosines.}
  	\label{coslaw}
  \end{figure}

 Consider now the general triangle in Figure \ref{coslaw},
  with the vectors $\ba $, $\bb $, $\bc $ along its sides satisfying the vector equation $ \ba + \bb = \bc $. Squaring this equation gives

 \[  (\ba + \bb )^2 =\ba^2 + \ba \bb + \bb \ba + \bb^2 =\bc^2\ \ \iff \ \ 
 |\ba|^2 + 2 \ba \cdot \bb + |\bb|^2 = |\bc|^2 , \]
 known as the {\it Law of Cosines}, where  
 \beq \ba \cdot \bb := \frac{1}{2}(\ba \bb + \bb \ba ) = |\ba ||\bb |\cos \theta ,\label{reverseab} \eeq
  is the {\it inner product} or {\it dot product} of the
 vectors $\ba $ and $\bb $. In Figure \ref{coslaw}, the angle $-\pi \le \theta \le \pi$ is measured from the vector
 $\ba $ to the vector $\bb $, and 
  \[\cos \theta = - \cos (\pi -\theta) =-\cos C = \cos (-\theta), \]
 so the sign of the angle is unimportant. Note that (\ref{reverseab}) allows us to reverse the order of the geometric product,
 \beq \bb \ba = - \ba \bb + 2 \ba \cdot \bb . \label{reverseba} \eeq   

We have used the usual rules for the multiplication of real
 numbers, except that we have not assumed that the multiplication of vectors is
 universally commutative. Indeed, the Pythagorean Theorem tells us that $|\ba|^2+ |\bb|^2 = |\bc|^2$ only for a right triangle when $\ba \cdot \bb =0$,
 or equivalently, when the vectors $\ba $ and $\bb $ are orthogonal and anti-commute.
 
 Now is a good place to summarize the rules which we have developed for the geometric multiplication of vectors. For vectors $\ba $, $\bb $, and $\bc$, 
 \renewcommand{\labelenumi}{(P\arabic{enumi})}
 \begin{enumerate}
 	\item $ \ba^2 = |\ba |^2 $ \quad {\it The square of a vector is its magnitude squared} \item $\ba \bb =- \bb \ba $ \quad {\it defines the bivector $\bB =\ba \bb $ when $\ba $ and $\bb $ are orthogonal vectors.}  
 	\item
 	$\ba ( \bb + \bc ) =  \ba\bb + \ba\bc$ \quad {\it Left distributivity}
   
 	\item
 	$( \bb + \bc)\ba = \bb\ba + \bc\ba$ \quad {\it Right distributivity}
 	\item     
 	$\ba (\bb\bc ) =(\ba\bb )\bc = \ba\bb\bc$     \quad {\it Product associativity}
 	\item $0 \ba  = 0 = \ba 0 $   \quad {\it Multiplication of a vector by zero is zero}
 	\item
 	$\alpha\ba  =  \ba\alpha, \ \  {\rm for} \ \ \alpha \in \R$  \quad {\it Multiplication of a vector times a scalar is commutative}
  \end{enumerate}
 
\subsection{The outer product} 
 So far, all is well, fine and good. The inner product of two vectors has been identified
 as one half the symmetric product of those vectors. To discover the geometric interpretation of the anti-symmetric product of the two vectors $\ba $ and $\bb $, we write
 \beq \ba \bb = \frac{1}{2}(\ba \bb + \bb \ba)+ \frac{1}{2}(\ba \bb - \bb \ba) =
    \ba \cdot \bb + \ba \w \bb , \label{geoproductab} \eeq
  where $\ba \w \bb :=  \frac{1}{2}(\ba \bb - \bb \ba)$ is called
  the {\it outer product}, or {\it wedge product} between $\ba $ and $\bb $. The outer
  product is {\it antisymmetric}, since $\bb \w \ba = - \ba \w \bb$.
  Indeed, when $\ba \cdot \bb =0$ the geometric product reduces to the outer product, {\it i.e.} 
  \beq \ba \bb = \ba \cdot \bb + \ba \w \bb = \ba \w \bb = - \bb \ba . \label{reducedprod} \eeq 
  
  It is natural to give $\ba \w \bb$ the interpretation of a {\it directed plane segment} or {\it bivector}. To see this, write $\bb = \bb_\parallel + \bb_\perp$,
  where $\bb_\parallel = s \ba $, for $s\in \R$, is the vector part of $\bb$ which is parallel to $\ba$, and $\bb_\perp$ is the vector part of $\bb $ which is perpendicular to $\ba$. Calculating $\ba \bb $, we find
  \[ \ba \bb = \ba (\bb_\parallel + \bb_\perp   ) = \ba \bb_\parallel   +\ba  \bb_\perp   =
   s \ba ^2 + \ba \bb_\perp  = \ba \cdot \bb + \ba \w \bb. \]
   Equating scalar and bivector parts, gives
   \beq \ba \cdot \bb = s \ba^2 \quad {\rm and} \quad \ba \w \bb = \ba \bb_\perp.  
                    \label{perpbivector} \eeq
   It follows that $\ba \w \bb = \ba \bb_\perp $ is the bivector which is the product of the
   orthogonal vectors $\ba $ and $\bb_\perp$, shown in Figure \ref{coslaw}.
      The bivector defined by the oriented parallelogram $\ba \w \bb $, with sides $\ba $ and $\bb $, has exactly the same orientation and directed area as the bivector  defined by the oriented rectangle $\ba \bb_\perp$, with the sides
   $\ba $ and $\bb_\perp$.  
     	\begin{figure}[h]
     		\centering
     		\includegraphics[width=0.7\linewidth]{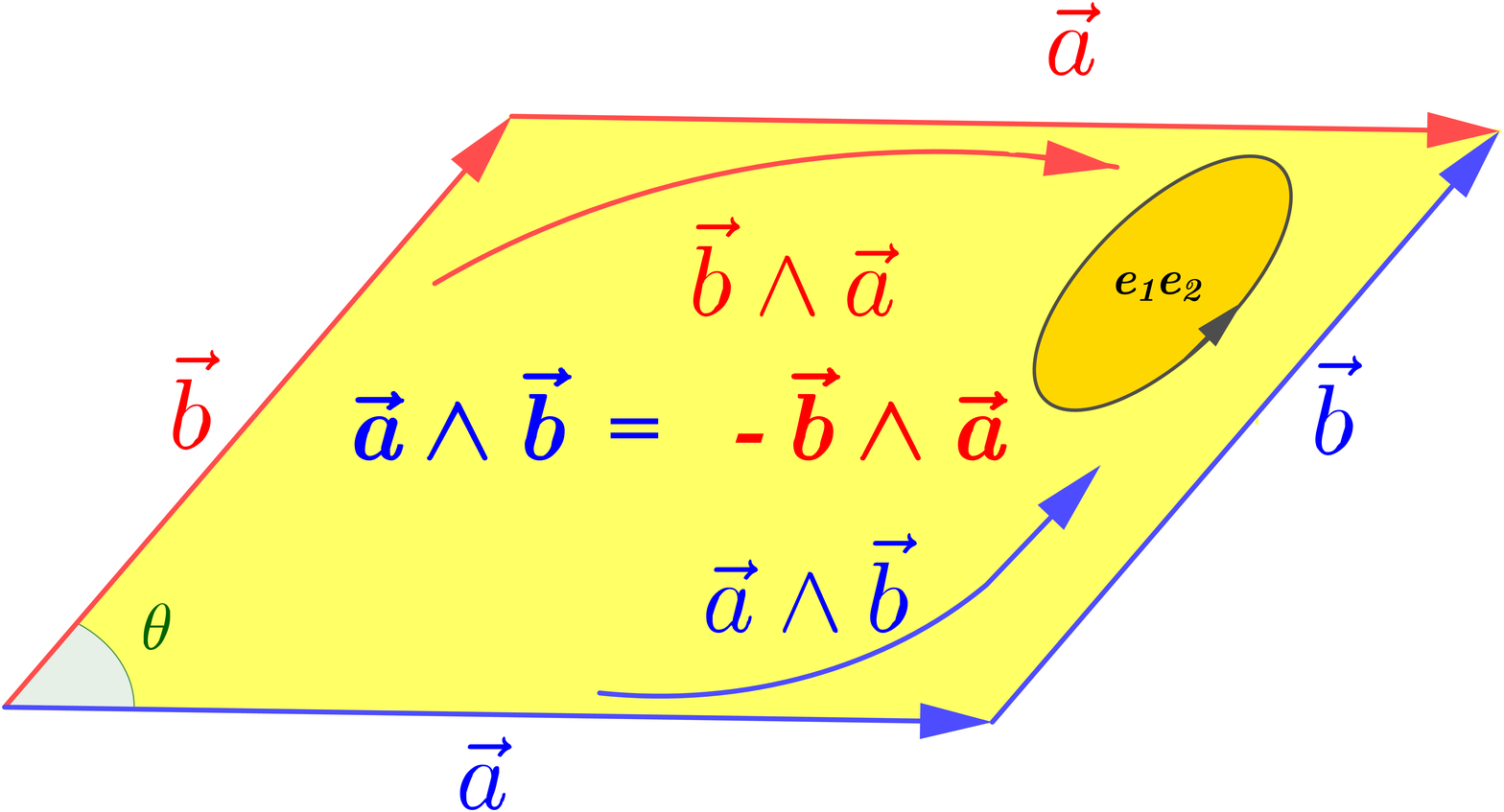}
     		\caption[Figure 3]{Orientation of a bivector. The {\it area}, or {\it magnitude}  of the bivector $\ba \w \bb$ is $|\ba \w \bb|= |\ba||\bb ||\sin \theta|$, where
     			$-\pi \le \theta < \pi$, and its direction is the unit bivector $\be_1 \be_2$. Note that the {\it shape} of the bivector $\be_{12}:=\be_1 \be_2$ is unimportant, only the plane in which it lies and its orientation.}
     	     		\label{bivecorient}
     	\end{figure}
    
      We have seen that
       the square of a vector is its magnitude squared, $\ba ^2 = |\ba |^2$.     
      What about the square of the bivector $(\ba \w \bb)$? Using (\ref{perpbivector}), we find that
      \beq (\ba \w \bb )^2 = (\ba \bb_\perp)^2 = - \ba^2 \bb_\perp^2   =-|\ba \w \bb|^2, \label{areabivector2}  \eeq
       in agreement with (\ref{bivectorsquared}). If the bivector is in the $xy$-plane of the unit bivector $\be_1  \be_2$, where the unit vectors $\be_1$ and $\be_2$ lie along the orthogonal $x$- and $y$-axes, respectively, then
       $\ba \w \bb =\be_{12}|\ba ||\bb | \sin \theta$, see Figure \ref{bivecorient}.  The geometric product in $\R^2$ and $\R^3$ is further discussed in Section 3.  
      
  \begin{figure}[h]
  	\centering
  	\includegraphics[width=0.7\linewidth]{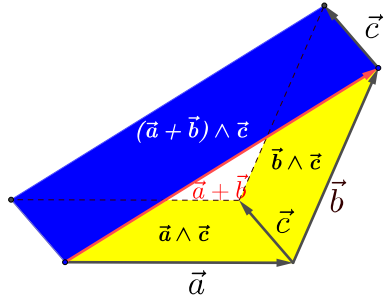}
  	\caption{The wedge product is distributive over the addition of vectors.}
  	\label{fig:distpropwedge}
  \end{figure}
   
   Just as sum of vectors is a vector,
   the sum of bivectors is a bivector. Figure \ref{fig:distpropwedge} shows 
   the sum of the bivectors
   \[ \ba \w \bc + \bb \w  \bc =( \ba + \bb )\w \bc, \]
   and also shows the {\it distributive property} of the outer product 
   over the sum of the
   vectors $\ba$ and $\bb$.

 \subsection{Properties of the inner and outer products}

 Since the triangle in Figure \ref{coslaw} satisfies the vector equation
 \[  \ba + \bb = \bc, \] 
 by wedging both sides of this equation by $\ba, \bb $ and $\bc$, gives
 \[ \ba \w \bb = \bc \w \bb, \ \ \bb \w \ba = \bc \w \ba , \ \ {\rm and} \ \ \bc \w \ba =
    \bb \w \bc , \]
    or equivalently,
  \[  \ba \w \bb = \bc \w \bb = \ba \w \bc . \]
   Note that the area of the triangle is given by $\frac{1}{2}|\ba \w \bb|$,
    which is one half of the area of the parallelogram $\ba \w \bb $, so the last equation is reflecting the equivalent relationship between parallelograms. 
     
\begin{figure}[h]
 	\centering
 	\includegraphics[width=0.55\linewidth]{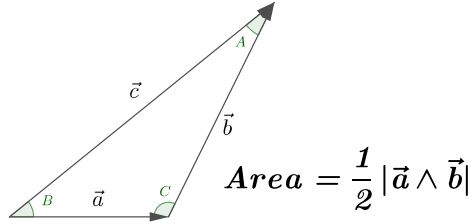}
 	\caption{Law of Sines.}
 	\label{sinlaw}
 \end{figure}

Dividing each term of the last equality by $|\ba || \bb ||\bc |$, gives
 \[ \frac{\hat \ba \w \hat \bb }{|\bc |} = \frac{\hat \bc \w \hat \bb }{|\ba |} = 
 \frac{\hat \ba \w \hat \bc }{|\bb |} \ \ \implies  \frac{|\hat \ba \w \hat \bb| }{|\bc |} = \frac{|\hat \bc \w \hat \bb| }{|\ba |} = 
 \frac{|\hat \ba \w \hat \bc |}{|\bb |}. \]
 For the angles $0 \le A,B,C \le \pi$,
 \[ |\hat \ba \w \hat \bb| =\sin C = \sin(\pi -C), \ \  |\hat \bc \w \hat \bb|= \sin A, \ \ {\rm and} \ \  |\hat \ba \w \hat \bc| = \sin B, \]
 from which it follows that
    \[  \frac{\sin A}{|\ba |} =\frac{\sin B}{|\bb| }=\frac{\sin C}{|\bc|} \]
   known as the {\it Law of Sines}, 
  see Figure \ref{sinlaw}.

 In (\ref{geoproductab}), we discovered that the geometric product of two vectors splits into two parts, a symmetric {\it scalar} part $\ba \cdot \bb $ and an anti-symmetric {\it bivector} part $\ba \w \bb$. It is natural to ask the question whether the geometric product of a vector $\ba$ with a bivector $\bb \w \bc $ has a similar decomposition? Analogous to (\ref{geoproductab}), we write
 \beq \ba (\bb \w \bc) = \ba \cdot (\bb \w \bc) + \ba \w (\bb \w \bc), \label{gaproductabc} \eeq
 where in this case
 \beq \ba \cdot (\bb \w \bc):=\frac{1}{2}\Big( \ba (\bb \w \bc) - (\bb \w \bc)\ba \Big)
 =: -  (\bb \w \bc)\cdot \ba 
  \label{adotbc} \eeq 
  is {\it antisymmetric}, and
  \beq \ba \w (\bb \w \bc):=\frac{1}{2}\Big( \ba (\bb \w \bc) + (\bb \w \bc)\ba \Big)
  =:    (\bb \w \bc)\w \ba 
    \label{awedgebc} \eeq 
  is {\it symmetric}.
      
\begin{figure}[h]
	\centering
		\includegraphics[width=0.7\linewidth]{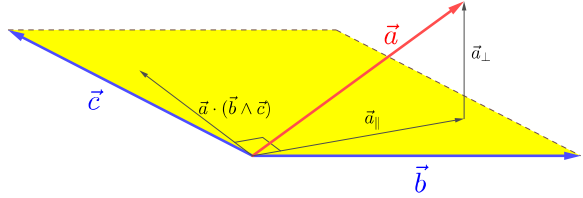}
	\caption{The result ${\ba}\cdot({\bb} \wedge {\bc})$ is the vector
		$\ba$ projected onto the plane of $\bb \w \bc $, and then rotated through $90$ degrees in this plane. }
	\label{apawc}
\end{figure}

  To better understand this decomposition, we consider each part separately. Starting with
  $\ba \cdot (\bb \w \bc)=\ba_\parallel (\bb \w \bc ) $, we first show that
  \beq \ba \cdot (\bb \w \bc) =(\ba \cdot \bb)\bc - (\ba \cdot \bc )\bb . \label{adotbc1} \eeq
  Decomposing the left side of this equation, using (\ref{adotbc}) and (\ref{geoproductab}), gives
    \[ \ba \cdot (\bb \w \bc) = \frac{1}{2}\Big(\ba (\bb \w \bc )-(\bb \w \bc )\ba \Big)
       = \frac{1}{4}(\ba \bb \bc - \ba \bc \bb -\bb  \bc \ba + \bc \bb \ba) . \]
 Decomposing the right side, gives
\[ (\ba \cdot \bb)\bc - (\ba \cdot \bc )\bb=\frac{1}{2}\Big( (\ba \cdot \bb)\bc + \bc (\ba \cdot \bb)     - (\ba \cdot \bc )\bb - \bb (\ba \cdot \bc )\Big)\]
\[=\frac{1}{4}\Big( (\ba  \bb + \bb \ba )\bc + \bc (\ba \bb + \bb \ba )     - (\ba  \bc+\bc \ba )\bb - \bb (\ba \bc + \bc \ba )\Big) \] 
\[= \frac{1}{4}(\ba \bb \bc - \ba \bc \bb -\bb  \bc \ba + \bc \bb \ba) ,    \]
which is in agreement with the left side. The geometric interpretation of (\ref{adotbc}) is given in the Figure \ref{apawc}.

\begin{figure}[h]
 	\centering
 	\includegraphics[width=0.7\linewidth]{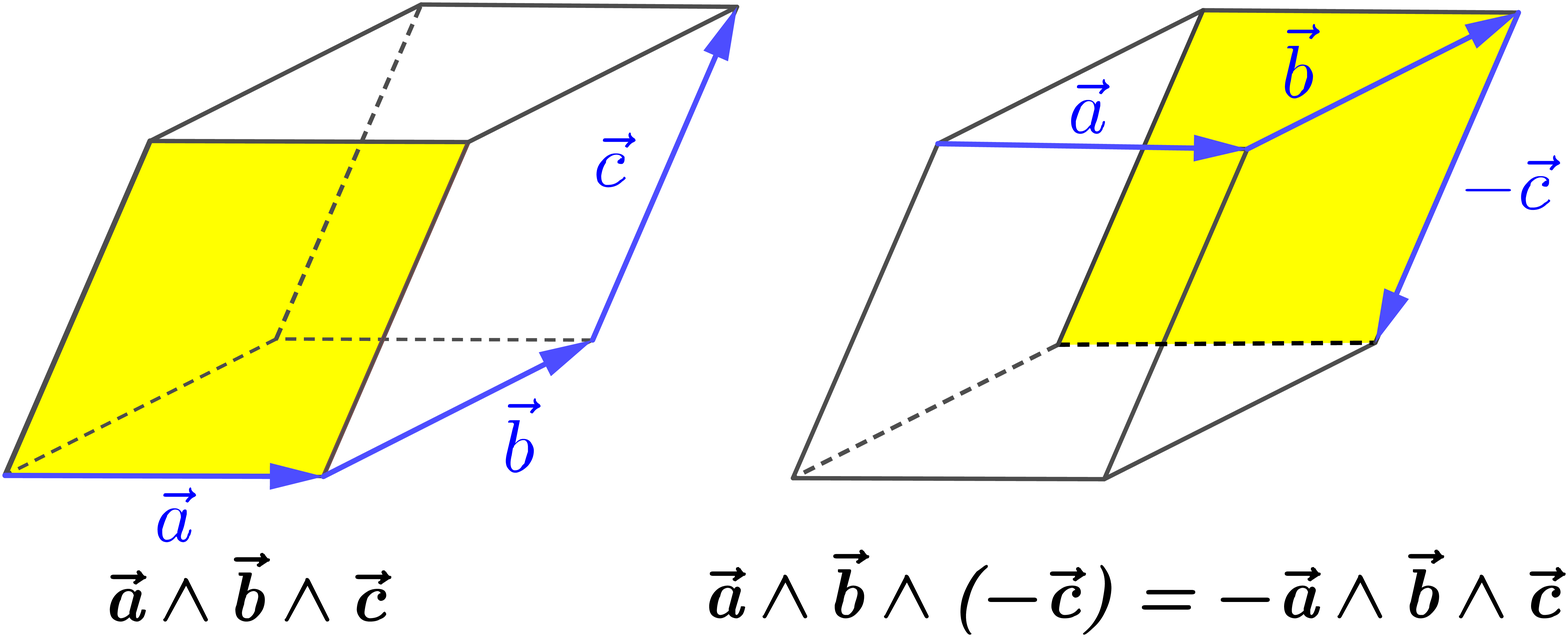}
 	\caption[Trivector]{The sign of the vector $\bc $ determines the {\it right}  and {\it left handed} orientation of the trivector $\ba \w \bb \w \bc$ shown.}
 	\label{trivector}
 \end{figure}

 Regarding the triple wedge product (\ref{awedgebc}), we need to show the associative
 	property, $\ba \w (\bb \w \bc)= (\ba \w \bb)\w \bc$. Decomposing both sides of this equation, using (\ref{trivector}) and (\ref{awedgebc}), gives
 \[ \ba \w (\bb \w \bc):=\frac{1}{2}\Big( \ba (\bb \w \bc) + (\bb \w \bc)\ba \Big)
 =  \frac{1}{4}\Big( \ba (\bb \bc- \bc \bb) + (\bb \bc- \bc \bb)\ba \Big) , \]
 and
 \[ (\ba \w \bb) \w \bc:=\frac{1}{2}\Big( (\ba \w \bb ) \bc + \bc (\ba \w \bb) \Big)
 =  \frac{1}{4}\Big(( \ba \bb - \bb \ba)\bc + \bc (\ba \bb - \bb \ba ) \Big) . \]
 To finish the argument, we have
 \[\ba \w (\bb \w \bc)- (\ba \w \bb)\w \bc = \frac{1}{4}\Big( - \ba \bc \bb -\bc \ba \bb +
 \bb \bc \ba +\bb\ba \bc \Big)  \]
 \[ =\frac{1}{2}\Big(-(\ba \cdot \bc )\bb + \bb (\ba \cdot \bc) \Big) =0. \]
 The {\it trivector} or {\it directed volume} $\ba \w \bb \w \bc$ is pictured in 
 Figure \ref{trivector}. 
  There are many more similar identities in higher dimensional geometric algebras \cite{H/S,SNF}.  
 
  \bigskip
 
 \no {\bf Exercise:} Using the properties (\ref{awedgebc}) and (\ref{adotbc1}), prove the Associative Law (P5) for the geometric product of vectors,
 \[   \ba ( \bb \bc )= (\ba \bb )\bc.  \]
   
 \section{The geometric algebras $\G_1$, $\G_2$ and $\G_3$.}

 In the previous section, we discovered two general principals for the multiplication
 of Euclidean vectors $\ba $ and $\bb $: 
  \begin{itemize}
 	\item[1)] The square of a vector is its length
 squared, $\ba^2 = |\ba |^2$.  
  \item[2)] If the vectors $\ba $ and $\bb $ are orthogonal to each other, i.e., the angle between
 them is $90$ degrees, then they anti-commute $\ba \bb = -\bb \ba $ and define the
 bivector given in (\ref{defbivector}). 
  \end{itemize}
 These two general rules hold for Euclidean vectors, independent of the dimension      
 of the space in which they lie. 
 
 The simplest euclidean geometric algebra is obtained by extending the real number system $\R $ to include a single new square 
 root of $+1$, giving the geometric algebra
 \[  \G_1 := \R(\be),  \]
 where $\be^2=1$. A geometric number in $\G_1$ has the form
\[  g= x + y \be,  \]
where $x,y \in \R $, and defines the hyperbolic number plane \cite{S1}.
 
\begin{figure}[h]
 	\centering
 	\includegraphics[width=0.75\linewidth]{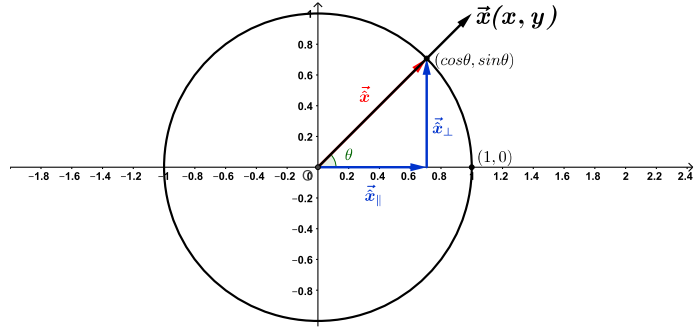}
 	\caption{The unit circle $S^1$ in the $xy$-plane.}
 	\label{unitcircle}
 \end{figure}

 We now apply what we have learned about the general geometric addition and multiplication of vectors to vectors in the two 
 dimensional plane $\R^2$, and in the three dimensional space $\R^3$ of experience. The $2$-dimensional {\it coordinate plane} is defined by 
 \beq \R^2 := \{(x,y)| \ \ x,y \in \R \}. \label{coorplane2} \eeq
 By laying out two {\it orthonormal unit vectors} $\{\be_1, \be_2\}$ along the $x$- and $y$-axes, respectively, each point 
 \beq (x,y)\in \R^2 \quad \longleftrightarrow \quad \bx = x \be_1 + y \be_2\in \R^2 \label{abusexy} \eeq 
 becomes a {\it position vector} $\bx = |\bx |\hat \bx  $ from the origin, shown in
 Figure \ref{unitcircle} with the unit circle. The point $\hat\bx  =(\cos \theta, \sin \theta)$ on the unit circle $S^1$, where the
 angle $\theta$ is measured from the $x$-axis, becomes the unit vector
 \[ \hat \bx=\cos (\theta)\be_1  +\sin (\theta)\be_2   . \]
 In equation (\ref{abusexy}), we have abused notation by equating the coordinate point $(x,y)\in \R^2$ with the {\it position vector} $\bx = x \be_1 + y \be_2$ from the origin of $\R^2$.      

\begin{figure}[h]
 	\centering
 	\includegraphics[width=0.7\linewidth]{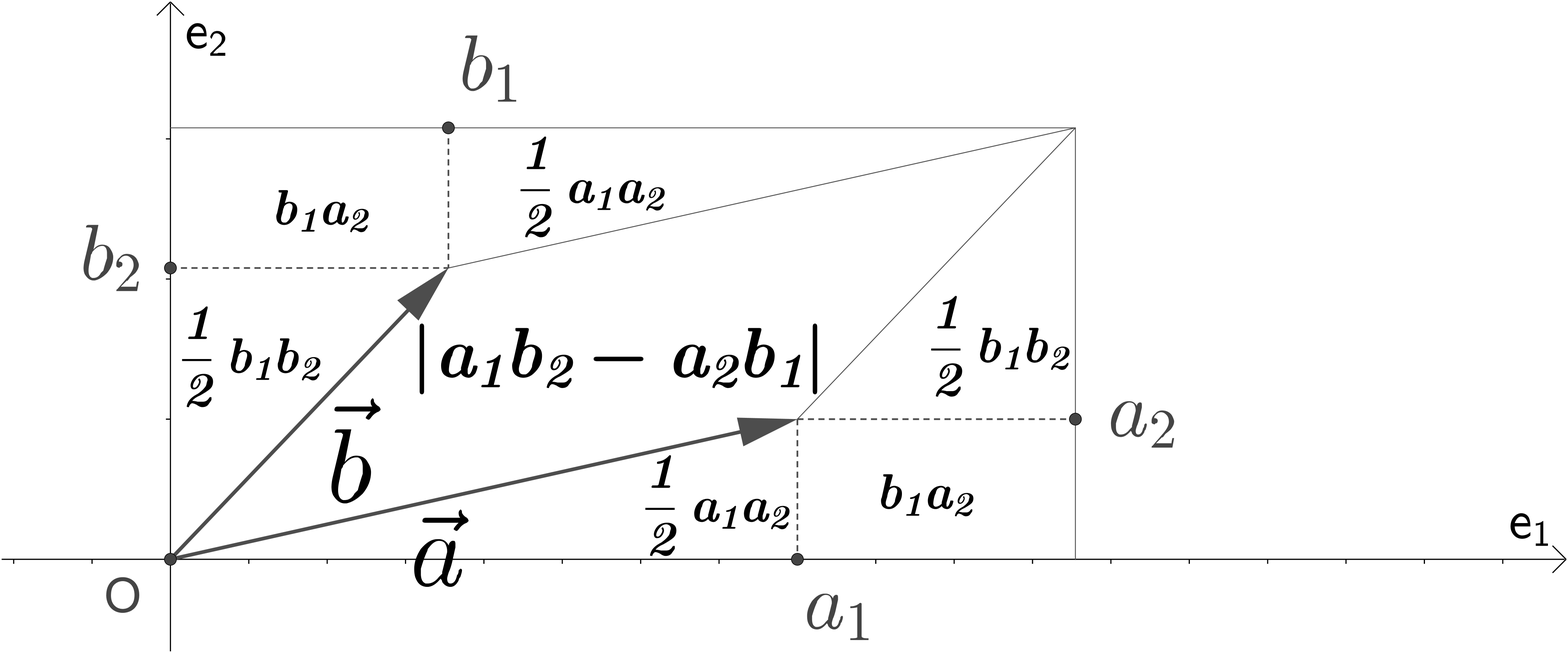}
 	\caption{The outer product $\ba \w \bb $ in 2-dimensions.}
 	\label{geoprod}
 \end{figure}

 Calculating the geometric product of the two vectors 
 \[ 	\ba=(a_1,a_2)=a_1 \be_{1} + a_2 \be_{2}, \quad 
 \bb=(b_1,b_2)=b_1 \be_{1} + b_2 \be_{2}, \]
 in the $xy$-plane, we obtain
 \[ \ba \bb = (a_1 \be_{1} + a_2 \be_{2})(b_1 \be_{1} + b_2 \be_{2}) \] 		
 \[ = a_1 b_1 \be_1^2 + a_2 b_2 \be_2^2 + a_1b_2 \be_1 \be_2+a_2 b_1\be_2 \be_1   \] 
 \beq = (a_1b_1 +a_2 b_2) +(a_1 b_2 -a_2 b_1 )\be_1 \be_2 = \ba \cdot \bb + \ba \w \bb  ,\label{dotouter2} \eeq 
 where the inner product $\ba \cdot \bb =  a_1b_1 +a_2 b_2 =|\ba||\bb|\cos \theta$, and the
 outer product
 \[ \ba \w \bb =(a_1 b_2 -a_2 b_1 )\be_{12} = \be_{12} |\ba ||\bb | \sin \theta \] 
 for $\be_{12}:= \be_1 \be_2 = \be_1 \wedge \be_2$.  
 The bivector $\ba \w \bb$ is pictured in Figure \ref{geoprod}, together
 with a picture proof that the magnitude $|\ba \w \bb|= |a_1b_2 -a_2 b_1|$, as expected.  
 
 By introducing the unit vectors $\{\be_1,\be_2\}$ along the coordinate axes of  $\R^2$, and using properties of the geometric product, we have found explicit formulas for the dot and outer products of any two vectors $\ba$ and $\bb$ in $\R^2$. The geometric product of the orthogonal unit vectors $\be_1$ and $\be_2$ gives the unit bivector 
 $\be_{12}$, already pictured in Figure \ref{bivecorient}. Squaring $\be_{12}$, gives
 \[ \be_{12}^2 = (\be_1 \be_2)(\be_1 \be_2)= -\be_1^2 \be_2^2 = -1,  \] 
 which because of (\ref{bivectorsquared}) and (\ref{areabivector2}) is no surprise.
 
 The most general geometric number of the $2$-dimensional Euclidean
 plane $\R^2$ is
 \[  g= g_0 + g_1 \be_1 + g_2 \be_2 + g_3 \be_{12}, \]  
 where $g_\mu \in \R$ for $\mu = 0,1,2,3$. The set of all geometric numbers $g$, together with the two operations of geometric addition and multiplication, make up the {\it geometric algebra} $\G_2$ of the Euclidean plane $\R^2$,
 \[  \G_2 :=    \{ g| \ \   g= g_0 + g_1 \be_1 + g_2 \be_2 + g_3 \be_{12}\} =\R(\be_1,\be_2)      . \]
 The formal rules for the geometric addition and multiplication of the geometric numbers in $\G_2$
 are exactly the same as the rules for addition and multiplication of real numbers, except we give up universal commutativity to express the
  anti-commutativity of orthogonal vectors.     
 
 The geometric algebra $\G_2$ breaks into two parts,
 \[ \G_2 = \G_2^0 + \G_2^1 + \G_2^2 =\G_2^+ + \G_2^-, \] 
 where the {\it even part}, consisting of {\it scalars} (real numbers) and bivectors,  
 \[ \G_2^+ :=\G_2^{0+2}=  \{ x + y\be_{12}| \ \ x,y \in \R \}\ \widetilde{=} \ \C  \]      
 is algebraically closed and isomorphic to the complex number $\C$, 
 and the {\it odd part}, 
 \[ \G_2^- := \G_2^1 = \{ \bx | \ \ \bx = x \be_1 + y \be_2    \}\equiv \R^2     \]
 for $x,y \in \R$, consists of vectors in the $xy$-plane $\R^2$.  
 The geometric algebra $\G_2$ unites the vector plane $\G_2^-$  and the complex number  plane $\G_2^+$ into a unified geometric number system $\G_2$ of the plane. 
  
 By introducing a third unit vector $\be_3$ into $\R^2$, along
 the $z$-axis, we get the $3$-dimensional space $\R^3$. All of the formulas found in $\R^2$ can then
 be extended to $\R^3$, and by the same process, to any higher $n$-dimensional space $\R^n$ for $n>3$. Geometric algebras can always be extended to higher dimensional geometric algebras simply by introducing additional orthogonal anti-commuting unit vectors with square $\pm 1$, \cite{GeoReal,Hyprevisit}.

\begin{figure}[h]
 	\centering
 	\includegraphics[width=0.7\linewidth]{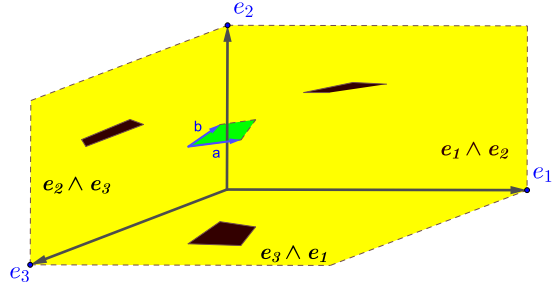}
 	\caption{Bivector decomposition in 3D space.}
 	\label{3dbivector}
 \end{figure}

Let us see how the
 formulas (\ref{dotouter2}) work out explicitly in  
 \beq  \R^3 := \{ \bx| \quad \bx = (x,y,z ) = x \be_1+y \be_2 + z \be_3 \}, \label{defR3} \eeq
 for $x,y,z \in \R$. For vectors   
 \[ \ba = a_1 \be_1 + a_2 \be_2+a_3 \be_3, \quad \bb = b_1 \be_1 + b_2 \be_2+b_3 \be_3,\]
we calculate
 \[  \ba \bb =( a_1 \be_1 + a_2 \be_2+a_3 \be_3) ( b_1 \be_1 + b_2 \be_2+b_3 \be_3) \] 
 \[ = a_1b_1 \be_1^2 + a_2 b_2 \be_2^2 + a_3 b_3 \be_3^2   \]
 \[ + a_1b_2 \be_1 \be_2 + a_2 b_1 \be_2 \be_1   + a_2 b_3 \be_2 \be_3  + a_3 b_2 \be_3 \be_2 +  a_1 b_3 \be_1 \be_3  + a_3 b_1 \be_3 \be_1  \]
 \[ =( a_1b_1 + a_2 b_2  + a_3 b_3) + ( a_1b_2  - a_2 b_1) \be_{12}   +
 ( a_2 b_3   - a_3 b_2) \be_{23} +  (a_1 b_3  - a_3 b_1) \be_{13}      \] 
 where the {\it dot} or {\it inner product}, 
 \[ \ba \cdot \bb =a_1b_1 + a_2 b_2  + a_3 b_3=|\ba ||\bb |\cos \theta , \] 
 and the outer product (\ref{perpbivector}), 
 \beq \ba \w \bb =( a_1b_2  - a_2 b_1) \be_{12}   +
 ( a_2 b_3   - a_3 b_2) \be_{23} +  (a_1 b_3  - a_3 b_1) \be_{13}=|\ba ||\bb |\hat \bB \sin \theta .  \label{awbarea} \eeq
 The sum of the three bivector components, which are projections onto the coordinate planes, are shown in Figure \ref{3dbivector}.
 
  In $\R^3$, the outer product $\ba \w \bb $ can be expressed in terms of the well known {\it cross product} of the century old, pre-Einstein Gibbs-Heaviside vector analysis.
 The vector cross product of the vectors $\ba $ and $\bb$ is defined by
 \[   \ba \times \bb := \det \pmatrix{\be_1 & \be_2 & \be_3 \cr
 	a_1 & a_2 & a_3 \cr b_1 & b_2 & b_3  }
 = ( a_2b_3  - a_3 b_2) \be_{1}   -
 ( a_1 b_3   - a_3 b_1) \be_{2} +  (a_1 b_2  - a_2 b_1) \be_{3} \] 
 \beq = |\ba || \bb | \sin \theta \, \hat \bn ,\label{vec-cross-product} \eeq
where $ \hat \bn :=
 \frac{\ba \times \bb }{|\ba \times \bb |}$ . 

\begin{figure}[h]
 	\centering
 	\includegraphics[width=0.7\linewidth]{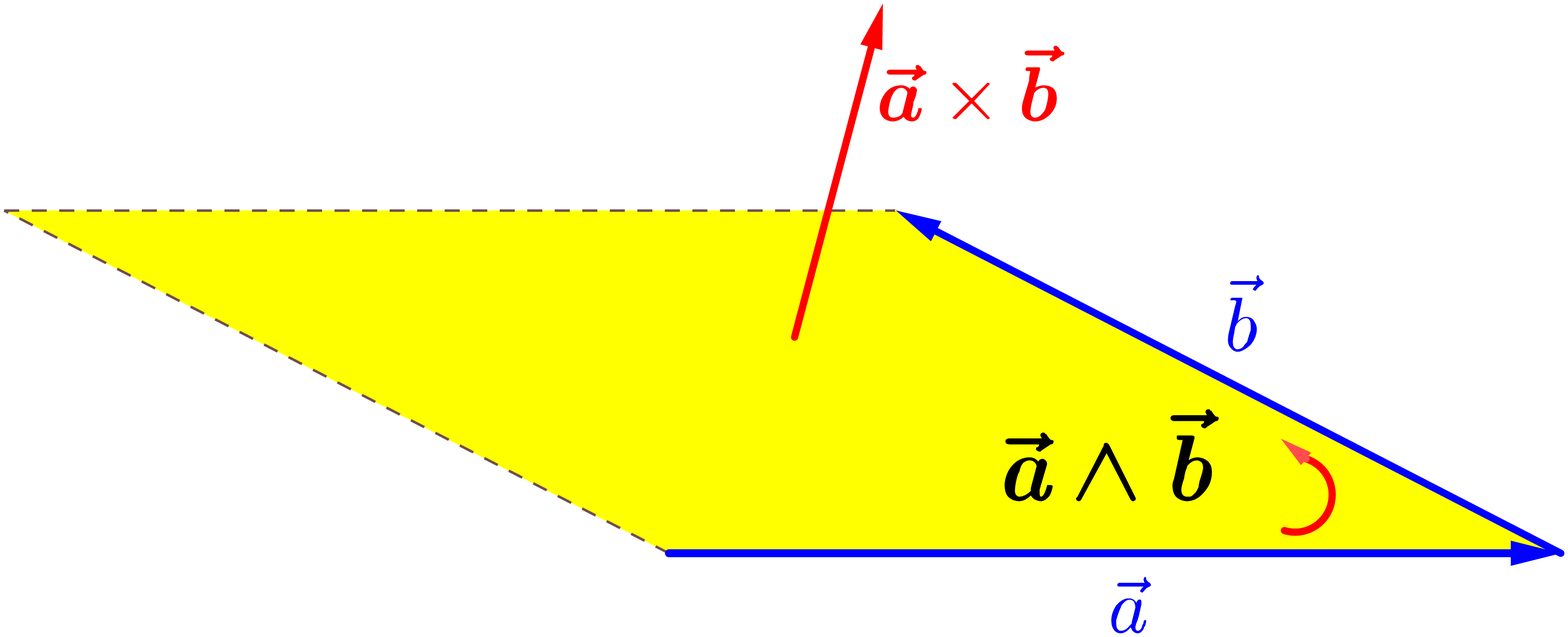}
 	\caption{The vector cross product $\ba \times \bb$ is the {\it right hand normal} dual to the
 		bivector $\ba \w \bb =I(\ba \times \bb) $. Also, $|\ba \times \bb |=|\ba \w \bb |$.}
 	\label{crossproductab}
 \end{figure}

 Defining the {\it unit trivector} or {\it pseudoscalar} of $\R^3$, 
  \beq  I := \be_1  \be_2 \be_3 = \be_{123},   \label{pseudoi} \eeq
 the relationship (\ref{vec-cross-product}) and (\ref{pseudoi}) can be combined into 
 \beq \ba \w \bb = I (\ba \times \bb)  = |\ba || \bb | \sin \theta \, I \hat \bn , \label{crosswedge} \eeq
  as can be easily verified. We say that
 the vector $\ba \times \bb$ is {\it dual} to, or the {\it right hand normal} of, the bivector $\ba \w \bb$, shown in the Figure \ref{crossproductab}.	Note that we are using the symbol $I=\be_{123}$ for the unit trivector or {\it pseudoscalar} of $\G_3$ to
 distinguish it from the $i=\be_{12}$, the unit bivector of $\G_{2}$.         
 
 We have seen in (\ref{geoproductab}) that the geometric product of two vectors
 decomposes into two parts, a scalar part and a vector part. We now calculate the
 geometric product of three vectors $\ba,\bb, \bc$.
 \[ \ba \bb \bc = \ba (\bb \cdot \bc + \bb \w \bc ) = (\bb \cdot \bc) \ba + \ba (\bb \w \bc) \]
 \[ =(\bb \cdot \bc) \ba + \ba\cdot (\bb \w \bc) + \ba \w \bb \w \bc . \]
 This shows that geometric product of three vectors consists of a vector part
 \[ (\bb \cdot \bc )\ba + \ba \cdot (\bb \w \bc) = (\bb \cdot \bc )\ba + (\ba \cdot \bb) \bc - (\ba \cdot \bc )\bb  \]
 and the trivector part $\ba \w \bb \w \bc$. For the vectors
 \[  \ba = a_1 \be_1 + a_2 \be_2 + a_3 \be_3, \ \bb = b_1 \be_1 + b_2 \be_2 + b_3 \be_3, \ \bc = c_1 \be_1 + c_2 \be_2 + c_3 \be_3 ,  \]
 the trivector 
 \beq  \ba \w \bb \w \bc = \det{\pmatrix{a_1 & a_2 & a_3 \cr
 		b_1 & b_2 & b_3 \cr c_1 & c_2 & c_3}} \be_{123}
 =\Big(( \ba \times \bb) \cdot \bc \Big)I  . \label{detabc} \eeq  
 
 By the {\it standard basis} of the
 geometric algebra $\G_3$ of the $3$-dimensional Euclidean space $\R^3$, we mean
 \[  \G_3 := span_\R \{ 1, \be_1 ,\be_2,  \be_3, \be_{12},\be_{13}, \be_{23}, \be_{123}     \}= \R(\be_1, \be_2, \be_3)  . \]  
 A general geometric number of $\G_3$ is 
 \[  g= g_0 + \bv + B + T \]
 where $g_0 \in \R $, $\bv = v_1 \be_1 + v_2 \be_2 + v_3 \be_3$ is a vector, 
 $B = b_{12}\be_{12} + b_{23}\be_{23} + b_{13}\be_{13}$ is a bivector,
 and $T=t I$, for $t\in \R $, is a {\it trivector} or {\it directed volume element}. Note that just like
 the unit bivector $i=\be_{12}$ has square $i^2=-1$, the unit trivector
 $I=\be_{123}$ of space has square $I^2=-1$, as follow from the calculation
 \[ I^2 = (\be_1 \be_2 \be_3)(\be_1 \be_2 \be_3) =(\be_1 \be_2) (\be_1 \be_2) 
  \be_3^2 = (-1)(+1)=-1.     \] 
  Another important property of the pseudoscalar $I$ is that it commutes with all vectors in $\R^3$, and hence with all geometric numbers in $\G_3$.  
 
\section{Analytic Geometry}

\begin{figure}[h]
	\centering
	\includegraphics[width=0.55\linewidth]{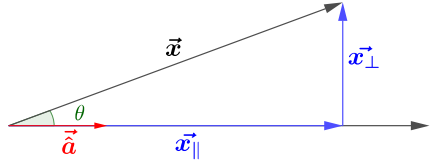}
	\caption{The vector $\bx $ is decomposed into parallel and perpendicular commponents with respect to the vector $\hat \ba $.} %\textcolor{red}{Remove the vector $\vec{a}$ from the figure.}}
	\label{decompvec}
\end{figure}
 
Given a vector $\bx$ and a unit vector $\hat \ba$, we wish to express $\bx =\bx_\parallel + \bx_\perp $ where $\bx_\parallel$ is parallel to $\hat \ba$, and $\bx_\perp$ is perpendicular
to $\hat \ba $, as shown in Figure \ref{decompvec}. 
Since $\hat \ba \hat \ba= 1$, and using the associative law, 
	\beq
	\nonumber
 \bx = (\bx \hat \ba) \hat \ba = (\bx \cdot \hat \ba )\hat \ba + (\bx \w \hat \ba )\hat \ba 
= \bx_\parallel + \bx_\perp, \label{vecdeomp}	\eeq
where 
 \[ \bx_\parallel = (\bx \cdot \hat \ba )\hat \ba \quad {\rm and} \quad \bx_\perp = (\bx \w \hat \ba)\hat \ba = \bx - \bx_\parallel . \]
We could also accomplish this decomposition by writing
\[  \bx =\hat \ba ( \hat \ba \bx)  =  \hat \ba (\hat \ba \cdot \bx )  + \hat \ba (\hat \ba \w \bx ) 
= \bx_\parallel + \bx_\perp . \]
It follows that $\bx_\parallel = (\bx \cdot \hat \ba) \hat \ba =\hat \ba (\bx \cdot \hat \ba )  $ as expected, and
 \[ \bx_\perp = (\bx \w \hat \ba )\cdot \hat \ba = \hat \ba \cdot (\hat \ba \w \bx )= 
     -\hat \ba \cdot (\bx \w \hat \ba)            , \]
 in agreement with (\ref{adotbc1}).

\begin{figure}[h]
	\centering
	\includegraphics[width=0.7\linewidth]{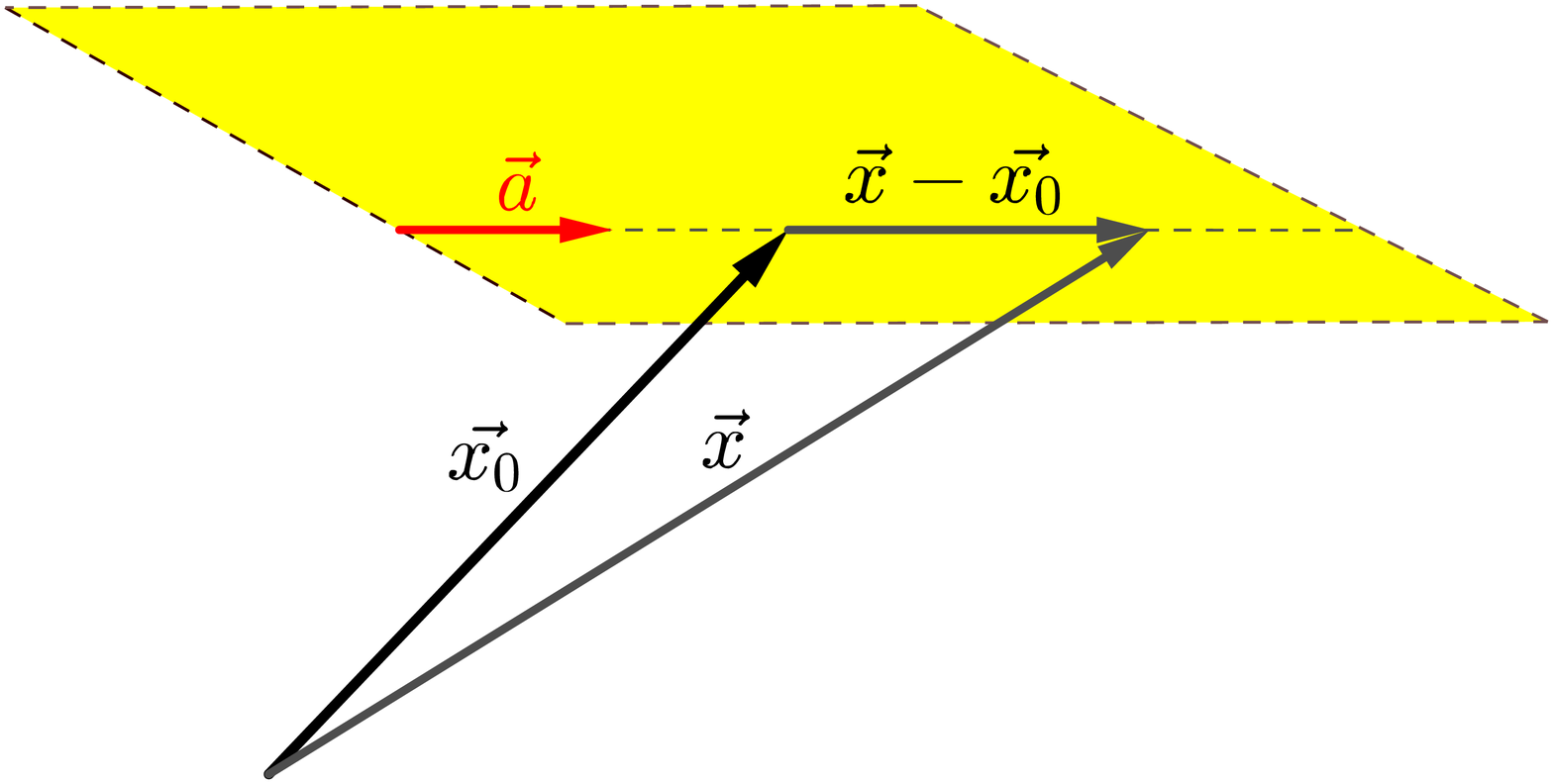}
	\caption{The line $L_{\bx_0}(\ba)$ through the point $\bx_0$ in the direction $\ba$.}
	\label{lineeq}
\end{figure}

One of the simplest problems in analytic geometry is given a vector $\ba $ and a point $\bx_0$, what is the equation of the line passing through the point $\bx_0$ in the
direction of the vector $\ba$? The line $L_{\bx_0}(\ba)$ is given by
\[  L_{\bx_0}(\ba):= \{ \bx| \ \ (\bx-\bx_0)\w \ba =0  \} .\]
The equation 
\[ (\bx-\bx_0)\w \ba =0 \ \ \iff \ \ \bx = \bx_0 + t \ba ,   \]
for $t \in \R$, see Figure \ref{lineeq}. 
\begin{figure}[h]
	\centering
	\includegraphics[width=0.6\linewidth]{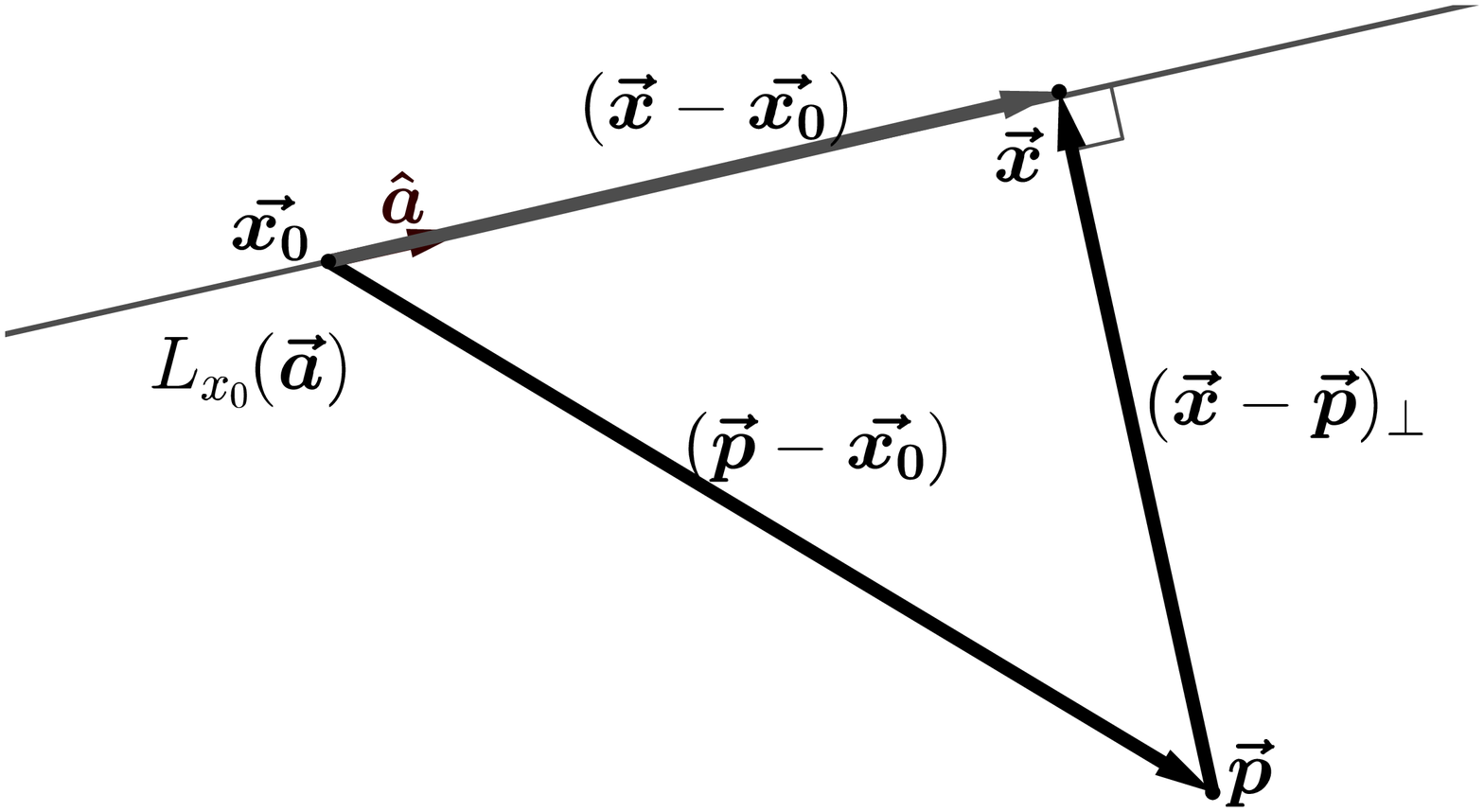}%{closestlinept}
	\caption{The distance of the point $\bp $ from the line $L_{\bx_0}(\ba)$ is 
		$| \bx - \bp |$.}
	\label{distance}
\end{figure}

Given the line $L_{\bx_0}(\ba )$, and a point $\bp $, let us find the point $\bx $ on the line $L_{\bx_0}(\ba )$ which is closest to the point $\bp$, and the distance $|\bx - \bp |$ from $\bx $ to $\bp $. Referring to Figure \ref{distance}, and using the decomposition (\ref{vecdeomp}) to project $\bp - \bx_0$ onto the vector
$\hat \ba $, we find
\[  \bx = \bx_0 + [(\bp - \bx_0 )\cdot \hat \ba  ] \hat \ba ,  \]
so, with the help of (\ref{geoproductab}) and (\ref{vecdeomp}),
\beq \bx - \bp=
 (\bx_0 - \bp )- [(\bx_0 - \bp )\cdot \hat \ba  ] \hat \ba  =(\bx_0 - \bp  )_\perp ,  \label{perpdist} \eeq
 where $(\bx_0 - \bp)_\perp$ is the component of $\bx_0 - \bp$ perpendicular to $\ba $.   
Using (\ref{perpdist}), the distance of the point $\bp$ to the line is
\[  | \bx - \bp|  =\sqrt{(\bx - \bp)^2} 
= \sqrt{(\bx_0-\bp )^2- \big((\bx_0 - \bp )\cdot \hat \ba\big)^2}=| (\bx_0 - \bp )_\perp|,  \]
  see Figure \ref{distance}.

\subsection{The exponential function and rotations}

The Euler exponential function arises naturally from the geometric product (\ref{geoproductab}). With the help of (\ref{reverseab}) and (\ref{crosswedge}), and noting that $(I\hat\bn)^2 = -1 $, the
geometric product of two unit vectors $\hat\ba$ and $\hat\bb $ in $\R^3$ is
\beq  \hat \ba \hat \bb = \hat \ba \cdot \hat \bb + \hat \ba \w \hat \bb =
\cos \theta +I \hat \bn  \sin \theta = e^{\theta I \hat \bn  }, \label{abrot} \eeq
where $\cos \theta := \hat \ba \cdot \hat \bb $. 
Similarly, 
\beq  \hat \bb \hat \ba = \hat \bb \cdot \hat \ba + \hat \bb \w \hat \ba =
\cos \theta -I \hat \bn  \sin \theta = e^{-\theta I \hat \bn  }. \label{barot} \eeq

\begin{figure}[h]
	\centering
	\includegraphics[width=0.65 \linewidth]{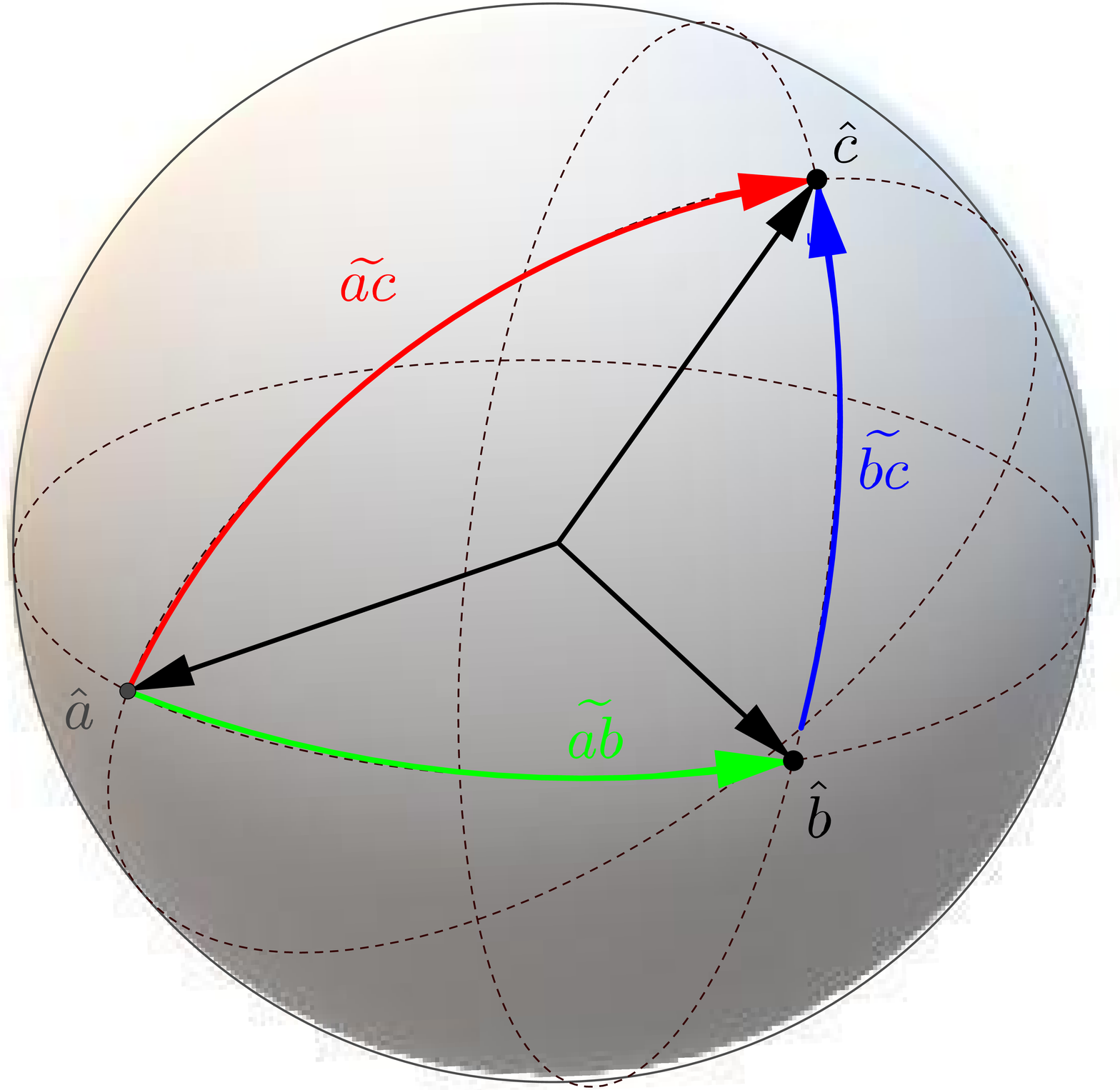}
	\caption{On the unit sphere, the arc $\widetilde{\hat\ba \hat\bb}$, followed by the arc $\widetilde{\hat\bb \hat\bc}$, gives the arc $\widetilde{\hat \ba \hat\bc}$. }
	\label{rota}
\end{figure}

Let $\hat \ba , \hat \bb , \hat \bc $ be unit vectors in $\R^3$. The equation
\beq (\hat \bb \hat \ba ) \hat \ba = \hat \bb (\hat \ba \hat \ba ) = \hat \bb =  (\ba \hat \ba) \hat \bb = \hat \ba (\hat \ba \hat \bb),   \label{fullangle} \eeq 
shows that when $\hat \ba $ is multiplied on the right by $\hat \ba \hat \bb = e^{\theta I \hat \bn}$, or on the left by $\hat \bb \hat \ba = e^{-\theta I \hat \bn}$,  it rotates the vector $\hat \ba $ through the angle $\theta$ into the vector $\hat \bb $.  
 The composition of rotations, can be pictured as the composition of arcs on the unit sphere. The composition of the arc $\widetilde{\bold{\hat{a}} \bold{\hat{b}}}$ on the great circle connecting the points $\hat\ba$ and $\hat\bb$, with the arc $\widetilde{\hat \bb \hat \bc }$ connecting $\hat\bb$ and $\hat\bc $, gives the arc $\widetilde{\hat \ba \hat \bc}$ connecting  $\hat\ba$ and $\hat\bc $. Symbolically,
\[ \widetilde{\hat \ba \hat\bb}\widetilde{\hat\bb \hat\bc} := (\hat\ba \hat \bb) (\hat \bb \hat \bc)=\hat \ba \hat \bc =: \widetilde{\hat\ba \hat\bc} ,  \] 
as shown in
\begin{figure}[h]
	\centering
	\includegraphics[width=0.7\linewidth]{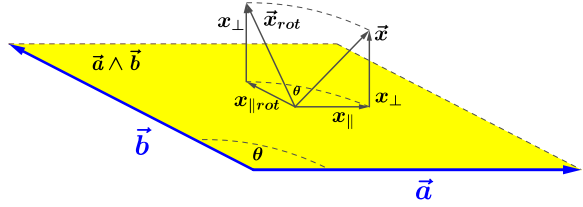}
	\caption{The parallel component $\bx_\parallel $ of $\bx $ in the plane of $\ba \w \bb $ is rotated through the angle $\theta$, leaving the perpendicular component $\bx_\perp$ unchanged.} 
	\label{xparallelrot}
\end{figure}  Figure \ref{rota}.  

By taking the square roots of both sides of equations (\ref{abrot}) and (\ref{barot}), it follows that 
\[ \sqrt{\hat \ba \hat \bb } = e^{\frac{1}{2}\theta I \hat \bn },  \quad {\rm and} \quad \sqrt{\hat \bb \hat \ba } = e^{- \frac{1}{2}\theta I \hat \bn }. \]   
Note also that
\beq \hat \bb = (\hat \bb \hat \ba) \hat \ba =
(\sqrt{\hat\bb \hat \ba })^2 \hat \ba =
\sqrt{\hat\bb \hat \ba }\ \hat \ba \, \sqrt{\hat\ba \hat \bb }
=  e^{-\frac{1}{2}\theta I \hat \bn  }\hat \ba  \,  e^{\frac{1}{2}\theta I \hat \bn  }.  \label{halfanglerot} \eeq
The advantage of the equation (\ref{halfanglerot}) over (\ref{fullangle})
is that it can be applied to rotate any vector $\bx$. For $\bx =\bx_\parallel + \bx_\perp $, where $\bx_\parallel$ is in the plane of $\ba \w \bb $, and $\bx_\perp$ is perpendicular to the plane, we get with the help of (\ref{adotbc}) and (\ref{awedgebc}),
\beq \bx_{rot} :=\sqrt{\hat\bb \hat \ba }\, \bx \,   \sqrt{\hat\ba \hat \bb }= e^{-\frac{1}{2}\theta I \hat \bn  }\, ( \bx_\parallel + \bx_\perp )  e^{\frac{1}{2}\theta I \hat \bn  }
=  e^{-\theta I \hat \bn  } \bx_\parallel + \bx_\perp , \label{rotatex} \eeq     
see Figure \ref{xparallelrot}. Formula (\ref{rotatex}) is known as the
{\it half angle} representation of a rotation \cite[p.55]{SNF}. A rotation can
also be expressed as the composition of two reflections.

 \subsection{Reflections}
 
 A bivector characterizes the direction of a plane. The equation of a plane passing through the origin in the direction of the bivector $\ba \w \bb$ is
 \beq  Plane_0(\ba \w \bb) = \{ \bx | \ \ \bx \w \ba \w \bb =0\}. \label{planeab} \eeq
 The condition that $\bx \w \ba \w \bb =0$ tells us that $\bx $ is in the the plane of the bivector $\ba \w \bb$, or 
 \[  \bx = t_a \ba + t_b \bb ,\]
 where $t_a, t_b \in \R$. This is the parametric equation of a plane passing through the origin having the direction of the bivector $\ba \w \bb$. If, instead, we want the
 equation of a plane passing through a given point $\bx_0$ and having the direction of
 the bivector $\ba \w \bb$, we have
 \beq  Plane_{\bx_0} (\ba \w \bb) = \{ \bx | \ \ (\bx-\bx_0) \w \ba \w \bb =0\}, \label{planeabx0} \eeq
 with the corresponding parametric equation
 \[ \bx = \bx_0 + t_a \ba + t_b \bb. \]

\begin{figure}[h]
   	\centering
   	\includegraphics[width=0.5\linewidth]{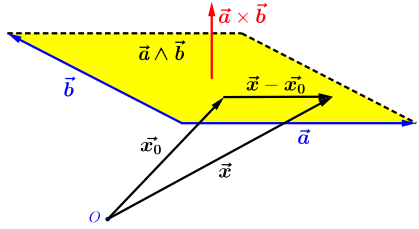}
   	\caption{The point $\bx$ is in the plane passing through the point $x_0$ and 
   		having the direction of the bivector $\ba \w \bb$.	} 
   	\label{planeeq}
   \end{figure}

  For a plane in $\R^3$, when $\bx=(x,y,z)$ and $\bx_0 = (x_0,y_0,z_0)$, using (\ref{detabc}) and (\ref{planeabx0}),
  \[  Plane_{\bx_0}(\ba \w \bb )
                   =  \{ \bx | \ \ \det{\pmatrix{x-x_0 & y-y_0 & z - z_0 \cr
                   		 a_1 & a_2 & a_3 \cr b_1 & b_2 & b_3} } =0\},  \]
   which is equivalent to the well known equation of a line through the point $\bx_0$,  
   \[ (\bx - \bx_0)\cdot \bn = 0 , \]
   where $\bn = \ba \times \bb $ is the {\it normal vector} to the bivector
   $\ba \w \bb $ of the plane, see Figure \ref{planeeq}. 
    
   Given a vector $\bx $ and a unit bivector $\ba \w \bb $, we decompose $\bx $ into a part $\bx_\parallel$ parallel to $\ba \w \bb $, and a part $\bx_\perp$   perpendicular to $\ba \w \bb$. Since by (\ref{awedgebc})
   \[ \bx_\parallel \w \ba \w \bb =\frac{1}{2}\Big(\bx_\parallel (\ba \w \bb)+(\ba \w \bb)\bx_\parallel \Big) = 0 ,\]
   and by (\ref{adotbc}),
    \[ \bx_\perp \cdot ( \ba \w \bb) =\frac{1}{2}\Big(\bx_\perp (\ba \w \bb)-(\ba \w \bb)\bx_\perp \Big) = 0, \]
  it follows that the parallel and perpendicular parts of $\bx $ anti-commute and
  commute, respectively, with the bivector $\ba\w \bb$. Remembering that
   $(\ba \w \bb )^2 = -1$, it follows that
  \beq  (\ba \w \bb )\bx (\ba \w \bb)= (\ba \w \bb )(\bx_\parallel + \bx_\perp)(\ba \w \bb)=\bx_\parallel - \bx_\perp \label{awb-mirror}. \eeq  
  This is the general formula for the reflection of a vector $\bx$ in a mirror
  in the plane of the unit bivector $\ba \w \bb$. 
  
  When we are in the $3$-dimensional space $\R^3$, the unit bivector 
  \[  \ba \w \bb = I (\ba \times \bb )= I \hat \bn . \]
   In this case, the reflection (\ref{awb-mirror}) takes the form
  \beq  (\ba \w \bb )\bx (\ba \w \bb)= - \hat \bn \bx \hat \bn = -\hat \bn (\bx_\parallel+ \bx_\perp)\hat \bn = \bx_\parallel - \bx_\perp. \label{axb-mirror} \eeq
 Since a rotation in $\R^3$ is generated by two consecutive reflections about two planes with normal unit vectors $\hat \bn_1 $ and $\hat \bn_2$, we have
  	\beq  \bx_{rot} = - \hat\bn_2 (- \hat \bn_1 \bx  \hat \bn_1) \hat \bn_2 = 
    (\hat \bn_2 \hat \bn_1 ) \bx (\hat \bn_1 \hat \bn_2 ).  \label{rotinR3} \eeq
%\textcolor{red}{\bf We need a new Figure ``xparallelrot" here.} 	

Letting $\hat\bn_1 \hat \bn_2=e^{\frac{1}{2} \theta I \hat \bn }$ 
where 
\[ \hat \bn := \frac{\hat \bn_1\times \hat \bn_2}{|\hat \bn_1\times \hat \bn_2| }, \] 
the formula for the rotation (\ref{rotinR3}) becomes
\beq  \bx_{rot} = (\hat \bn_2 \hat \bn_1 ) \bx (\hat \bn_1 \hat \bn_2 )
=e^{-\frac{1}{2} \theta I \hat \bn }\bx e^{\frac{1}{2} \theta I \hat \bn }
=e^{-\frac{1}{2} \theta I \hat \bn }\bx_\parallel e^{\frac{1}{2} \theta I \hat \bn } + \bx_\perp , \label{rotinR3a} \eeq
which is equivalent to (\ref{rotatex}).

 \section{Stereographic projection and a bit of quantum mechanics}
 \begin{figure}[h]
 	\begin{center}
 			\includegraphics[scale=.1]{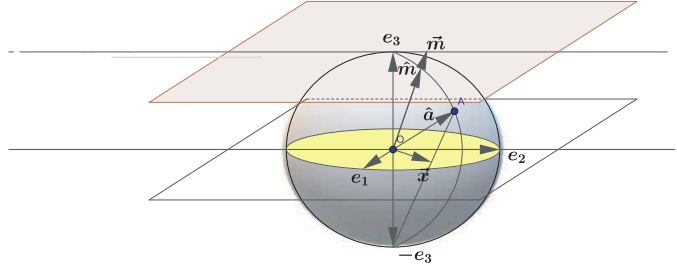}
 	 		\caption{Stereographic Projection from the South Pole of $S^2$ to the $xy$-plane,
 			where $\bm = \bx + \be_3$ and $\hat \bm = \frac{\bm}{|\bm |}$.}
 		\label{sterox}
 	\end{center}
 \end{figure}

\begin{figure}[h]
        	\begin{center}
        		\no\includegraphics[scale=.1]{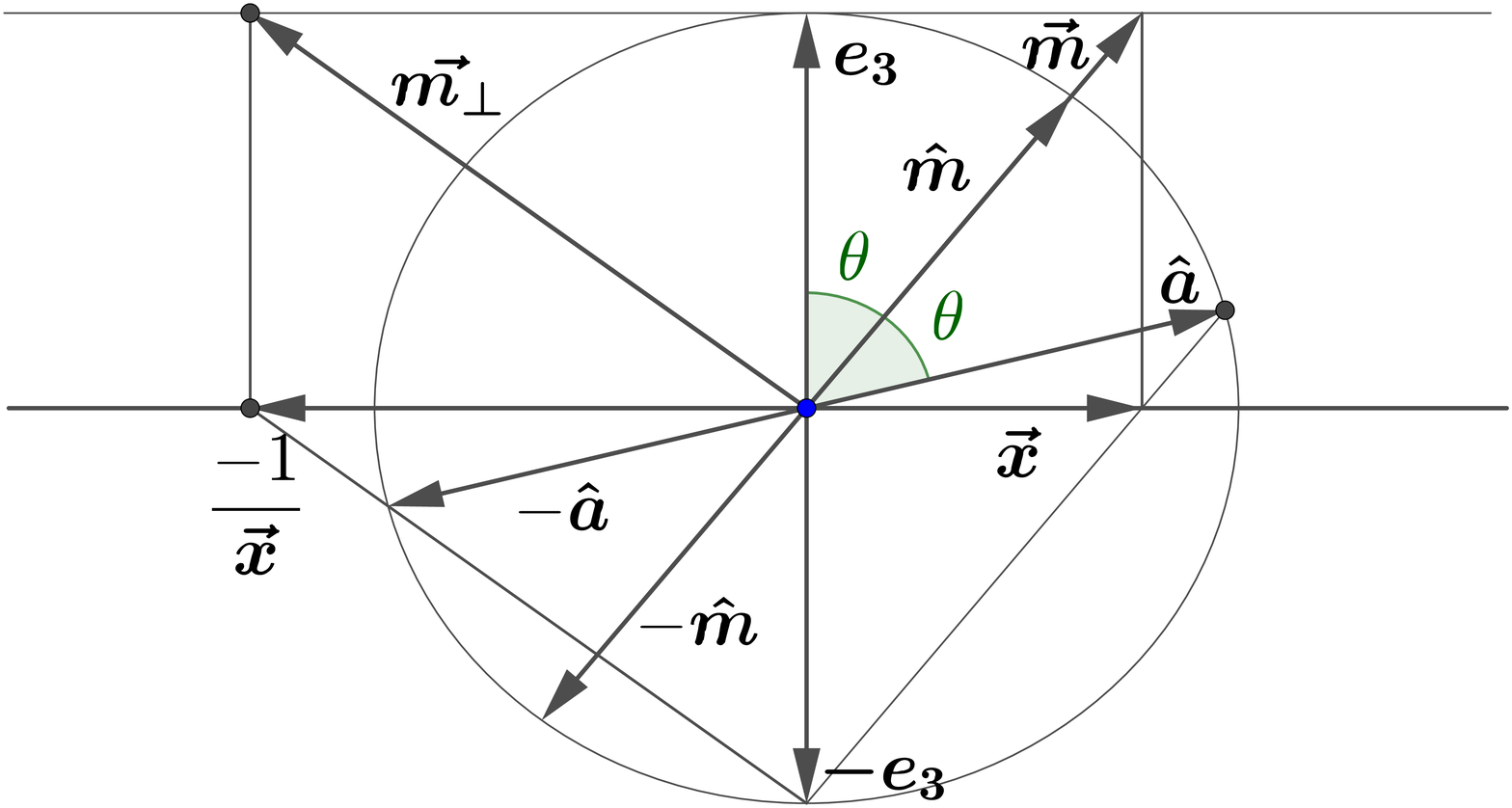}
        		\caption{A 2-D cut away in the plane of great circle through the points
        			$\be_3, \hat \ba $, and $-\hat \ba $ on $S^2$.  }
        		\label{sterox2}
        	\end{center}
        \end{figure}

As a final demonstration of the flexibility and power of geometric algebra, we  discuss stereographic projection from the unit sphere $S^2 \subset \R^3$, defined by 
\[ S^2 := \{\hat \ba | \quad \hat \ba^2 =1 \ \ {\rm and} \ \ \hat \ba \in \R^3  \}, \] 
onto $\R^2$. The mapping $\bx = f(\hat \ba)\in \R^2 $ defining stereographic projection is 
 \beq \bx = f(\hat \ba) := \frac{2}{\hat \ba + \be_3}-\be_3 , \ \ {\rm where} \ \
    \hat \ba \in S^2,\label{stereoproj}  \eeq
and is pictured in Figure \ref{sterox}. A 2-D cut away in the plane
of the great circle, defined by the points $\be_3, \hat \ba$, and the origin, is shown in Figure \ref{sterox2}.   
   Stereographic projection is an example of a {\it conformal mapping}, which preserves angles, and has many important applications in mathematics, physics, and more recently in robotics \cite{ECGS01,Sob2012}.  
        
 In working with the mapping (\ref{stereoproj}), it is convenient to use the new variable $\bm = \bx + \be_3$, for which case the mapping takes the simpler form
 \beq \bm = \frac{2}{\hat \ba + \be_3} = \frac{2(\hat \ba + \be_3)  }
 { (\hat \ba + \be_3)^2} = \frac{\hat \ba + \be_3  }{1+
 	\hat \ba \cdot \be_3}  .  \label{stereoprojm} \eeq 
 The effect of this change of variable maps points $\bx \in \R^3$ into corresponding
 points $\bm$ in the plane $Plane_{\be_3}(\be_{12})$ passing through the point $\be_3$ and parallel to $\R^2= Plane_0(\be_{12} )$. Noting that 
 \[ \be_3 \cdot \bm =\be_3 \cdot \Big( \frac{\hat \ba + \be_3  }{1+
 	\hat \ba \cdot \be_3}\Big)  = 1,     \]
 and solving the equation (\ref{stereoprojm}) for $\hat \ba $, gives with the help of (\ref{inversevec}) and (\ref{reverseba}),
 \[  \hat \ba = \frac{2}{\bm}-\be_3= \bm^{-1}\big(2- \bm \be_3 \big) \]
 \beq           = \frac{\hat \bm }{|\bm |}\big(2+\be_3 \bm -2 \be_3\cdot \bm \big)
 =\hat \bm \be_3 \hat \bm   . \label{ahateqn} \eeq
 We also have
 \beq  \hat \ba=\hat \bm \be_3 \hat \bm =(\hat \bm \be_3)\be_3(\be_3 \hat \bm) =(-I\hat \bm) \be_3 (I\hat \bm) , \label{ahateqn3} \eeq
 showing that $\hat \ba$ is obtained by a rotation of $\be_3$ in the plane of $\hat \bm\w \be_3$ through an
 angle of $2\theta$ where $\cos\theta:=\be_3\cdot \hat \bm$, or equivalently,
 by a rotation of $\be_3$ in the plane of $I \hat \bm$ through an angle of $\pi$.
  
  Quantum mechanics displays many surprising, amazing, and almost magical properties, which defy the classical mechanics of everyday experience.
 If the {\it quantum spin state} of an electron is put into a spin state $\hat \ba \in S^2 $ by a strong magnetic field at a given time, then the {\it probability of observing the electron's spin} in the spin state $\hat\bb \in S^2$ at a time immediately thereafter is
 \beq prob_{\hat \ba}^+(\hat \bb ) := \frac{1}{2}(1 + \hat \ba \cdot \hat \bb )=
 1- \frac{(\bm_a - \bm_b)^2 }{\bm_a^2 \bm_b^2 }, \label{probab}  \eeq
 where 
 \[ \hat \ba = \frac{2}{\bm_a }- \be_3 \quad {\rm and} \quad  
  \hat \bb = \frac{2}{\bm_b }- \be_3,  \]
  see \cite{SNov16,SMar17}. 
  
\begin{figure}[b]
  	\begin{center}
  		\no\includegraphics[scale=0.41]{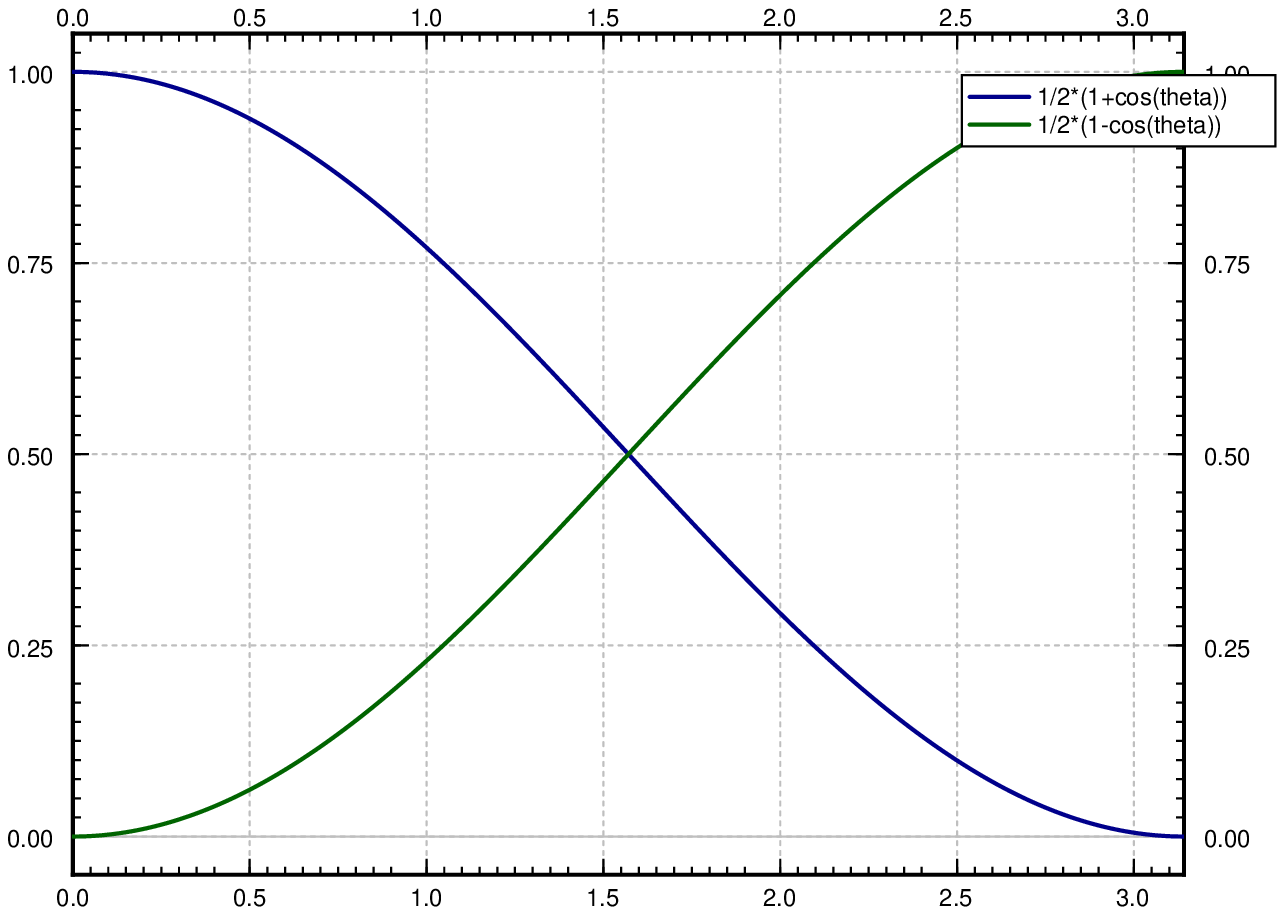}
  		\caption{The functions $prob_{\hat\ba }^\pm(\hat \bb )$. The angle 
  			$0 \le \theta \le \pi$ is between the unit vectors $\hat \ba $ and $\hat \bb $.} 
  		\label{plotplusminus}
  	\end{center}
  \end{figure}

  On the other hand, the {\it probability of a photon being emitted} 
  by an electron prepared in a spin state $\hat \bb $, when it is forced by a magnetic field into the spin state $\hat \ba $ is
    \beq prob_{\hat \ba}^-(\hat \bb ) := \frac{1}{2}(1 - \hat \ba \cdot \hat \bb )=
    \frac{(\bm_a - \bm_b)^2 }{\bm_a^2 \bm_b^2 }. \label{probabminus}  \eeq
 Whenever a photon is emitted, it has {\it exactly the same energy}, regardless of the angle
 $\theta$ between the spin states $\hat\ba$ and $\hat \bb $, \cite{SLN06,SLNY}.
 A plot of these two probability functions is given in Figure \ref{plotplusminus}.
 The equalities in (\ref{probab}) and (\ref{probabminus}) show that $prob_{\hat \ba }^\pm(\hat \bb )$ is directly related to the Euclidean distances between the points $\bm_a,\bm_b \in Plane_{\be_3}(\be_{12})$.
 %Clearly $prob_{\hat \ba }(\hat \ba)=1$, $prob_{\hat \ba }(\hat %\ba_\perp)=\frac{1}{2}$, and $prob_{\hat \ba }(-\hat \ba)=0$.
  The case when 
 \[ \bm_{a}=\bm=\bx + \be_3, \ \ \hat \bb =-\hat \ba, \ \  {\rm  and} \  \
  \bm_{b} = \bm_\perp = - \frac{1}{\bx}+ \be_3 \] 
  is pictured in Figure \ref{sterox2}.   

\end{document}

%% file: macros.tex
\def\ba{\mathbf{a}}   %\def\ba{{\bf a}}
\def\bb{\mathbf{b}}   %\def\bb{{\bf b}}
\def\bc{\mathbf{c}}   %\def\bc{{\bf c}}
\def\be{\mathbf{e}}   %\def\bd{{\bf e}}
\def\bm{\mathbf{m}}   %\def\bm{{\bf m}}
\def\bn{\mathbf{n}}   %\def\bn{{\bf n}}
\def\bp{\mathbf{p}}   %\def\bp{{\bf p}}
\def\bu{\mathbf{u}}   %\def\bu{{\bf u}}
\def\bv{\mathbf{v}}   %\def\bv{{\bf v}}
\def\bw{\mathbf{w}}   %\def\bw{{\bf w}}
\def\bx{\mathbf{x}}   %\def\bx{{\bf x}}
\def\bB{\mathbf{B}}   %\def\bB{{\bf B}}

\def\C{\mathbb{C}} %\def\C{I\!\!\!C} %\def\C{{C\kern-.647em I}}
\def\G{\mathbb{G}} % \def\G{I\!\!\!G}

\def\R{\mathbb{R}}  %\def\R{I\!\!R}

\def\no{\noindent}

\def\beq{\begin{equation}}
\def\eeq{\end{equation}}

\def\w{\wedge}
\def\implies{\Longrightarrow}

%% file: GA3mex__revised_.bbl
\begin{thebibliography}{}

  
  	\bibitem{WKC1882} W.K. Clifford, {\em Applications of Grassmann's extensive algebra}, Am. J. Math (ed.), Mathematical Papers by William Kingdon Clifford, pp. 397-401, Macmillan, London (1882). (Reprinted by Chelsea, New York, 1968.)
  	\bibitem{ECGS01} E.B. Corrochano, G. Sobczyk, Editors, {\em Geometric Algebra with Applications in Science and Engineering}, Birkh\"auser (2001).	\bibitem{MJC1985} M.J. Crowe, {\em A History of Vector Analysis}, Dover, New York (1985). 
  	\bibitem{TD1967} T. Dantzig, {\em Number: The Language of Science}, 4th edn. Free Press, New York 1967.
  	\bibitem{DL07} C. Doran, A. Lasenby, {\em Geometric Algebra for Physicists}, Cambridge 2007. 	
  	\bibitem{HD02} T. F. Havel, J.L. Doran, {\em Geometric Algebra in Quantum Information Processing},
  	Contemporary Mathematics, ISBN-10: 0-8218-2140-7, Vol. 305, 2002.
  	\bibitem{H/S} D. Hestenes and G. Sobczyk. {\it Clifford Algebra to
  		Geometric Calculus: A Unified Language for Mathematics and Physics},
  	2nd edition, Kluwer 1992.
  	\bibitem{LP97} P. Lounesto, 
  	\newblock {\em Clifford Algebras and Spinors, 2nd Edition}.
  	\newblock {Cambridge University Press}, Cambridge, 2001.
  	\bibitem{S1} G. Sobczyk, {\it Hyperbolic Number Plane}, The College Mathematics
  		Journal, Vol. 26, No. 4, pp.268-280, September 1995.
 	\bibitem{Sob2012} G. Sobczyk, {\it Conformal Mappings in Geometric Algebra}, Notices of the AMS, 
 		Volume 59, Number 2, p.264-273, 2012.
 	\bibitem{SNF} G. Sobczyk, {\em New Foundations in Mathematics: The Geometric Concept of Number},  
 	\newblock Birkh\"auser, New York 2013.
 	\bibitem{Shopf2015} G. Sobczyk, {\it Geometric Spinors, Relativity and the Hopf Fibration}, http://www.garretstar.com/geo2hopf26-9-2015.pdf
  
    \bibitem{SNov16} G. Sobczyk, {\it Part I: The Vector Analysis of Spinors} (2016)    
    http://www.garretstar.com/paulispin-19-07-2015.pdf \\
    https://arxiv.org/pdf/1507.06608.pdf
   	\bibitem{GeoReal} G. Sobczyk, {\it Geometrization of the Real Number System}, July 2017. http://www.garretstar.com/geonum2017.pdf
  		\bibitem{Hyprevisit} G. Sobczyk, {\it Hyperbolic Numbers Revisted}, Dec. 2017.
  		http://www.garretstar.com/hyprevisited12-17-2017.pdf
     \bibitem{SMar17} G. Sobczyk, {\it Spinors in Spacetime Algebra and Euclidean 4-Space} (2017) \\
     https://arxiv.org/pdf/1703.01244.pdf   
     \bibitem{SLN06} L. Susskind, {\it Lecture Notes 2:  Electron Spin}, Stanford University (2006). \\
      http://www.lecture-notes.co.uk/susskind/quantum-entanglements/lecture-2/electron-spin/
      \bibitem{SLNY} L. Susskind, {\it YouTube: Quantum Entanglements, Lecture 2} \\
      https://www.youtube.com/playlist?list=PL8D12D5AADF422C5D\&feature=plcp
 
  \end{thebibliography}
